\documentclass[12pt]{article}
\usepackage{mathrsfs}
\usepackage{amsmath,amssymb}

\usepackage{amsfonts}
\usepackage{latexsym}
\usepackage{amsthm}
\usepackage{pictex}

\usepackage[T2A]{fontenc}     
\usepackage[cp1251]{inputenc} 
\usepackage{graphicx}

\usepackage{epsfig}
\usepackage{fancyhdr}
\usepackage{color}

\numberwithin{equation}{section}

\newcommand{\er}[1]{{\rm(\ref{#1})}}
\def\lb{\label}

\newtheorem{theorem}{Theorem}[section]
\newtheorem{definition}{Definition}
\newtheorem{lemma}{Lemma}[section]

\newtheorem{proposition}{Proposition}[section]

\begin{document}

\def\a{\alpha} \def\cA{{\cal A}} \def\bA{{\bf A}}  \def\mA{{\mathscr A}}
\def\b{\beta}  \def\cB{{\cal B}} \def\bB{{\bf B}}  \def\mB{{\mathscr B}}
\def\g{\gamma} \def\cC{{\cal C}} \def\bC{{\bf C}}  \def\mC{{\mathscr C}}
\def\G{\Gamma} \def\cD{{\cal D}} \def\bD{{\bf D}}  \def\mD{{\mathscr D}}
\def\d{\delta} \def\cE{{\cal E}} \def\bE{{\bf E}}  \def\mE{{\mathscr E}}
\def\D{\Delta} \def\cF{{\cal F}} \def\bF{{\bf F}}  \def\mF{{\mathscr F}}
\def\c{\chi}   \def\cG{{\cal G}} \def\bG{{\bf G}}  \def\mG{{\mathscr G}}
\def\z{\zeta}  \def\cH{{\cal H}} \def\bH{{\bf H}}  \def\mH{{\mathscr H}}
\def\e{\eta}   \def\cI{{\cal I}} \def\bI{{\bf I}}  \def\mI{{\mathscr I}}
\def\p{\psi}   \def\cJ{{\cal J}} \def\bJ{{\bf J}}  \def\mJ{{\mathscr J}}
\def\vT{\Theta}\def\cK{{\cal K}} \def\bK{{\bf K}}  \def\mK{{\mathscr K}}
\def\k{\kappa} \def\cL{{\cal L}} \def\bL{{\bf L}}  \def\mL{{\mathscr L}}
\def\l{\lambda}\def\cM{{\cal M}} \def\bM{{\bf M}}  \def\mM{{\mathscr M}}
\def\L{\Lambda}\def\cN{{\cal N}} \def\bN{{\bf N}}  \def\mN{{\mathscr N}}
\def\m{\mu}    \def\cO{{\cal O}} \def\bO{{\bf O}}  \def\mO{{\mathscr O}}
\def\n{\nu}    \def\cP{{\cal P}} \def\bP{{\bf P}}  \def\mP{{\mathscr P}}
\def\r{\rho}   \def\cQ{{\cal Q}} \def\bQ{{\bf Q}}  \def\mQ{{\mathscr Q}}
\def\s{\sigma} \def\cR{{\cal R}} \def\bR{{\bf R}}  \def\mR{{\mathscr R}}
\def\S{\Sigma} \def\cS{{\cal S}} \def\bS{{\bf S}}  \def\mS{{\mathscr S}}
\def\t{\tau}   \def\cT{{\cal T}} \def\bT{{\bf T}}  \def\mT{{\mathscr T}}
\def\f{\phi}   \def\cU{{\cal U}} \def\bU{{\bf U}}  \def\mU{{\mathscr U}}
\def\F{\Phi}   \def\cV{{\cal V}} \def\bV{{\bf V}}  \def\mV{{\mathscr V}}
\def\P{\Psi}   \def\cW{{\cal W}} \def\bW{{\bf W}}  \def\mW{{\mathscr W}}
\def\o{\omega} \def\cX{{\cal X}} \def\bX{{\bf X}}  \def\mX{{\mathscr X}}
\def\x{\xi}    \def\cY{{\cal Y}} \def\bY{{\bf Y}}  \def\mY{{\mathscr Y}}
\def\X{\Xi}    \def\cZ{{\cal Z}} \def\bZ{{\bf Z}}  \def\mZ{{\mathscr Z}}
\def\O{\Omega}
\def\ve{\varepsilon}
\def\vt{\vartheta}
\def\vp{\varphi}
\def\vk{\varkappa}

\def\ba {\breve{a}}
\def\bb {\breve{b}}

\newcommand{\gA}{\mathfrak{A}} \newcommand{\ga}{\mathfrak{a}}
\newcommand{\gB}{\mathfrak{B}} \newcommand{\gb}{\mathfrak{b}}
\newcommand{\gC}{\mathfrak{C}}
\newcommand{\gD}{\mathfrak{D}}
\newcommand{\gE}{\mathfrak{E}}
\newcommand{\gF}{\mathfrak{F}}
\newcommand{\gG}{\mathfrak{G}}
\newcommand{\gH}{\mathfrak{H}}
\newcommand{\gI}{\mathfrak{I}}
\newcommand{\gJ}{\mathfrak{J}}
\newcommand{\gK}{\mathfrak{K}}
\newcommand{\gL}{\mathfrak{L}}
\newcommand{\gM}{\mathfrak{M}}
\newcommand{\gN}{\mathfrak{N}}
\newcommand{\gO}{\mathfrak{O}}
\newcommand{\gP}{\mathfrak{P}}
\newcommand{\gR}{\mathfrak{R}}  \newcommand{\gr}{\mathfrak{r}}
\newcommand{\gS}{\mathfrak{S}}   \newcommand{\gs}{\mathfrak{s}}
\newcommand{\gT}{\mathfrak{T}}
\newcommand{\gU}{\mathfrak{U}}
\newcommand{\gV}{\mathfrak{V}}
\newcommand{\gW}{\mathfrak{W}}
\newcommand{\gX}{\mathfrak{X}}
\newcommand{\gY}{\mathfrak{Y}}
\newcommand{\gZ}{\mathfrak{Z}}

\def\mA{{\mathscr A}}
\def\mB{{\mathscr B}}
\def\mC{{\mathscr C}}
\def\mD{{\mathscr D}}
\def\mE{{\mathscr E}}
\def\mF{{\mathscr F}}
\def\mG{{\mathscr G}}
\def\mH{{\mathscr H}}
\def\mI{{\mathscr I}}
\def\mJ{{\mathscr J}}
\def\mK{{\mathscr K}}
\def\mL{{\mathscr L}}
\def\mM{{\mathscr M}}
\def\mN{{\mathscr N}}
\def\mO{{\mathscr O}}
\def\mP{{\mathscr P}}
\def\mQ{{\mathscr Q}}
\def\mR{{\mathscr R}}
\def\mS{{\mathscr S}}
\def\mT{{\mathscr T}}
\def\mU{{\mathscr U}}
\def\mV{{\mathscr V}}
\def\mW{{\mathscr W}}
\def\mX{{\mathscr X}}
\def\mY{{\mathscr Y}}
\def\mZ{{\mathscr Z}}

\def\Z{{\Bbb Z}}
\def\R{{\Bbb R}}
\def\C{{\Bbb C}}
\def\T{{\Bbb T}}
\def\N{{\Bbb N}}
\def\S{{\Bbb S}}
\def\H{{\Bbb H}}
\def\J{{\Bbb J}}
\def\dD{{\Bbb D}}
\def\W{{\Bbb W}}
\def\K{{\Bbb K}}
\def\qqq{\qquad}
\def\qq{\quad}
\newcommand{\ma}{\begin{pmatrix}}
\newcommand{\am}{\end{pmatrix}}
\newcommand{\ca}{\begin{cases}}
\newcommand{\ac}{\end{cases}}
\let\ge\geqslant
\let\le\leqslant
\let\geq\geqslant
\let\leq\leqslant
\def\ma{\left(\begin{array}{cc}}
\def\am{\end{array}\right)}
\def\iint{\int\!\!\!\int}
\def\lt{\biggl}
\def\rt{\biggr}
\let\geq\geqslant
\let\leq\leqslant
\def\[{\begin{equation}}
\def\]{\end{equation}}
\def\wh{\widehat}
\def\wt{\widetilde}
\def\pa{\partial}
\def\sm{\setminus}
\def\es{\emptyset}
\def\no{\noindent}
\def\ol{\overline}
\def\iy{\infty}
\def\ev{\equiv}
\def\/{\over}
\def\ts{\times}
\def\os{\oplus}
\def\ss{\subset}
\def\h{\hat}
\def\Re{\mathop{\rm Re}\nolimits}
\def\Im{\mathop{\rm Im}\nolimits}
\def\supp{\mathop{\rm supp}\nolimits}
\def\sign{\mathop{\rm sign}\nolimits}
\def\Ran{\mathop{\rm Ran}\nolimits}
\def\Ker{\mathop{\rm Ker}\nolimits}
\def\Tr{\mathop{\rm Tr}\nolimits}
\def\const{\mathop{\rm const}\nolimits}
\def\dist{\mathop{\rm dist}\nolimits}
\def\diag{\mathop{\rm diag}\nolimits}
\def\cZr{\mathop{\rm Wr}\nolimits}
\def\BBox{\hspace{1mm}\vrule height6pt width5.5pt depth0pt \hspace{6pt}}

\def\Diag{\mathop{\rm Diag}\nolimits}


\def\Twelve{
\font\Tenmsa=msam10 scaled 1200 \font\Sevenmsa=msam7 scaled 1200
\font\Fivemsa=msam5 scaled 1200 \textfont\msbfam=\Tenmsb
\scriptfont\msbfam=\Sevenmsb \scriptscriptfont\msbfam=\Fivemsb

\font\Teneufm=eufm10 scaled 1200 \font\Seveneufm=eufm7 scaled 1200
\font\Fiveeufm=eufm5 scaled 1200
\textfont\eufmfam=\Teneufm \scriptfont\eufmfam=\Seveneufm
\scriptscriptfont\eufmfam=\Fiveeufm}

\def\Ten{
\textfont\msafam=\tenmsa \scriptfont\msafam=\sevenmsa
\scriptscriptfont\msafam=\fivemsa

\textfont\msbfam=\tenmsb \scriptfont\msbfam=\sevenmsb
\scriptscriptfont\msbfam=\fivemsb

\textfont\eufmfam=\teneufm \scriptfont\eufmfam=\seveneufm
\scriptscriptfont\eufmfam=\fiveeufm}

\title {Periodic Jacobi operator with finitely supported perturbation
on the half-lattice}

\author{ Alexei Iantchenko
\begin{footnote}
{ School of Technology,
  Malm{\"o} University,
SE-205 06 Malm{\"o}, Sweden, email: ai@mah.se }
\end{footnote} \and
Evgeny Korotyaev
\begin{footnote}
{Sankt Petersburg, e-mail: korotyaev@gmail.com}
\end{footnote}
}


\maketitle

\begin{abstract}
 \no We consider a periodic Jacobi operator $J$ with finitely supported perturbations  on the half-lattice.
 We describe all eigenvalues and resonances of $J$ and give
their properties. We solve the inverse resonance problem: we prove that the mapping from finitely supported
 perturbations to the Jost functions is one-to-one and onto, we show how the Jost functions can be reconstructed from all  eigenvalues, resonances and from the set of zeros of $S(\l)-1,$ where $S(\l)$ is the scattering matrix.
\end{abstract}

\section{Introduction.}
\setcounter{equation}{0}

We consider a  Jacobi operator $J=J^0+V$ on the half-lattice
$\N=\{1,2,3,..\}$. Here the unperturbed operator $J^0$ is a periodic
Jacobi operator given by
\[
\lb{J0}
 (J^0y)_n={a}_{n-1}^0y_{n-1}+{a}_{n}^0y_{n+1}+{b}_n^0y_n,\qqq n\ge1,\qqq
 y_0=0,
 \]
where  $y=(y_n)_1^\iy\in \ell^2=\ell^2(\N)$ and   the $q-$periodic
coefficients $a_n^0, b_n^0\in \R$ satisfy
  \[
 \lb{1e}
{a}^0_n={a}^0_{n+q}>0,\qq b_n^0={b}^0_{n+q},\qq n\in
\N=\{1,2,3....\},\qqq \prod_{j=1}^qa_j^0=1,\qq q\ge 2.
\]
The perturbation operator $V$ is the finitely supported Jacobi
operator given by
\[
\lb{V} (Vy)_n=\ca {u}_{n-1}y_{n-1}+u_{n}y_{n+1}+v_ny_n, & {\rm
if}\qq  1\le n\le p,\qq y_0=0,\\
u_py_p,& {\rm if} \qq n=p+1,\\
0, \qqq & {\rm if} \qq n\ge p+2,\qq p\ge 1.\ac
\]
We parameterize $V$ by the vector $(u,v)\in \R^{2p}$ and let $(u,v)$
belong to the class $\gX_\n$ given by
\begin{align}
\label{class}
 \gX_\nu=&\lt\{ (u,v)\in \R^{2p}:\,\,
a^0_n+u_n>0,\qq n=1,...,p, \ \ u_p\neq0\rt\}\qq \mbox{if}\qq \n=2p,
\\
\label{class2} \gX_\nu=&\lt\{ ({u} ,{v} )\in
\R^{2p}:\,\,{a}^0_n+{u} _n >0 ,\,\ \ n=1, ...,p,\ \ {v} _p\neq
 0,\,\,{u} _p= 0\rt\}\qq\mbox{if}\qq \n=2p-1.
\end{align}
We rewrite $J$ in the form
\[
\lb{pert}
 (J y)_n={a}_{n-1}y_{n-1}+{a}_{n}y_{n+1}+{b}_ny_n,\qqq n\ge1,\qqq
 y_0=0,
 \]
with the coefficients $a_n, b_n$ given by
\[
\lb{ab}
 {a}_n=\ca {a}^0_n+{u}_n>0& {\rm if} \qq n\le p,\\
                           {a}^0_n & {\rm if} \qq n\ge p+1, \ac
              \qqq {b}_n=\ca  b_n^0+{v}_n & {\rm if} \qq n\le p, \\
            b_n^0   & {\rm if} \qq n\ge p+1.         \ac
\]
 The  corresponding  Jacobi matrices have the forms
\[
\lb{Jm}
J^0=\left(\begin{array}{ccccc}
b_1^0   & a_1^0     &0         & 0        &... \\
  a^0_1 & b^0_2        &a^0_2& 0        &... \\
  0  &a^0_2& b^0_3         & a^0_3     &... \\
  0  & 0        & a^0_3      & b^0_4       &... \\
  0  & 0        &0         & a^0_4&... \\
  ... & ...      &...       &...       &... \\
 \end{array}\right),\qqq
   J=\left(\begin{array}{ccccc}
{b}_1   & a_1     &0         & 0        &... \\
a_1 & {b}_2        &{a}_2& 0        &... \\
  0  &{a}_2& {b}_3         & {a}_3     &... \\
  0  & 0        &{a}_3      & {b}_4       &... \\
  0  & 0        &0         &{a}_4&... \\
  ... & ...      &...       &...       &... \\
 \end{array}\right).
\]
 Note that  the $n=1$ case in (\ref{pert}) can be
thought of as forcing the Dirichlet condition $y_0=0.$ Thus,
eigenfunctions must be non-vanishing at $n=1$ and eigenvalues must
be simple.

The spectrum of $J^0$ consists of an
absolutely continuous part $\s_{\rm ac}(J^0)= \bigcup\limits_{1}^q
\s_j$ plus  at most one eigenvalue in each non-empty gap $\g_j,
j=1,...,q-1$, where the bands $\s_j$ and the gaps $\g_j$ are given
by
\begin{multline}
\lb{ends} \s_j=[\l_{j-1}^+,\l_j^-],\qq j=1,\ldots,q,\qq
\g_{j}=(\l^-_{j},\l^+_j),\qqq j=1,\ldots,q-1,\\
\l_0^+<\l_1^-\le \l_1^+<.. ...<\l^-_{q-1}\le \l^+_{q-1}
<\l^-_{q}.\qqq \qqq
\end{multline}
We introduce the infinite gaps
$$
\g_0=(-\infty, \l_0^+),\qqq \g_q=( \l_q^+,+\infty).
$$
Let $\vp=(\vp_n(\l))_1^\iy$ and $ \vt=(\vt_n(\l))_1^\iy$ be
fundamental solutions for the equation
\[
\lb{1e2} a_{n-1}^0 y_{n-1}+{a}_{n}^0y_{n+1}+{b}_n^0y_n=\l y_n,\qqq
\l\in\C,
\]
satisfying the conditions $\vt_0=\vp_1=1$ and $\vt_1=\vp_0=0$. Here and
below $a_0^0=a_q^0$.  Introduce the  Lyapunov function  $\D$ by
\[
\lb{LF} \D={\vp_{q+1}+\vt_{q}\/2}.
\]
It is known that $\D(\l)$ is a polynomial of degree $q$ and
$\l_j^\pm, j=1,...,q,$ are the zeros of the polynomial $\D^2(\l)-1$
of degree $2q$. Note that $\D(\l_j^\pm)=(-1)^{q-j}$. In each ``gap''
$[\l_j^-, \l_j^+]$ there is one simple zero of polynomials $\vp_q,
\dot \D, \vt_{q+1}.$ Here and below $\dot f$ denotes the derivative
of $f=f(\l)$ with respect to $\l:$ $\dot f\equiv\partial_\l f\equiv
f'(\l).$

Let $\G$ denote the complex plane cut along the segments $\sigma_j$ (\ref{ends}):
 $\G=\C\setminus\s_{\rm ac}(J^0).$ Now we introduce the two-sheeted Riemann
surface $\L$  of $\sqrt{1-\D^2(\l)}$  by joining the upper and lower rims of
two copies of the cut plane $\G$ in
the usual (crosswise) way.  We identify the first (physical) sheet $\L_1$ with $\G$
and the second sheet we denote by $\L_2$.

 Let $\,\,\wt{}\,\,$ denote  the natural projection from
 $\L$  into the complex plane $\C$:
 \begin{equation}
 \lb{projection}
\l\in\L,\qq\l\to \wt\l\in\C.
\end{equation}
By identification of $\G=\C\setminus\s_{\rm ac}(J^0)$ with $\L_1,$ the map $\,\,\wt{}\,\,$ can be also considered to be projection from $\L$ into the physical sheet $\L_1.$

The $j-$th gap on the first physical sheet $\Lambda_1$ we will denote by $\g_j^+$ and the same gap but on
the second nonphysical sheet $\Lambda_2$  we will denote by $\g_j^-$ and let $\g^{\rm c}_j$  be the union of $\overline{\g^+_j}$ and $\overline{\g^-_j}$:
\[
\lb{union}
\g_j^{\rm c}=\overline{\g^+_j}\cup\overline{\g^-_j}.
\]
Define the function $\O(\l)= \sqrt{1-\D^2(\l)}, \l\in \L,$ by
\begin{equation}\lb{branch}
\O(\l)<0\qq\mbox{for}\qq\l\in (\l_{q-1}^+, \l_{q}^-)\ss \L_1.
\end{equation} Introduce the
Bloch functions $\psi_n^\pm$ and the Titchmarch-Weyl functions
$m_\pm$ on $\L$ by
\begin{align}
 &\psi_n^\pm (\l)=\vt_n(\l)+m_\pm(\l)\vp_n(\l),\lb{Jost}\\
&m_\pm(\l)= \frac{\f(\l)\pm i\O(\l)}{\vp_q},\,\qq
\phi=\frac{\vp_{q+1}-\vt_q}{2},\,\,\l\in\L_1. \lb{TitchWeyl}
\end{align}
 The projection of all singularities of $m_\pm$ to the complex plane coincides with the set of zeros
$\{\mu_j\}_{j=1}^{q-1}$ of polynomial $\vp_q$. Recall that $\vt_n, \vp_n, \f$ are polynomials.
Recall that any polynomial  $P(\l)$  gives rise to a function $P(\l)=P(\wt\l)$  on
the  Riemann surface $\L$ of $\sqrt{1-\D^2(\l)}.$

The perturbation $V$ satisfying (\ref{V}) does not change the absolutely continuous spectrum:
\[
 \lb{acspectr}\s_{\rm ac}(J)=\s_{\rm ac}(J^0)=
 \bigcup_{j=1}^q[\l_{j-1}^+,\l_j^-].
\]
The spectrum of $J$ consists of an absolutely continuous part
$\s_{\rm ac}(J)=\s_{\rm ac}(J^0)$ plus a finite number of simple eigenvalues in
each non-empty gap $\g_j, j=0,...,q$.

In the present paper we consider the properties of the eigenvalues,
virtual states and resonances of the operators $J^0$ and $J,$ and solve
the inverse problem in terms of the resonances  of $J.$ Let
$R(\l)=(J-\l)^{-1}$ denote the resolvent of  $J$ and  let
$\langle\cdot,\cdot\rangle$ denote the scalar product in
$\ell^2=\ell^2(\N).$ Then for any $f,g\in\ell^2$ the function $\langle
R(\l)f,g\rangle$ is defined  on $\L_1$ outside the poles at the
bound states on the gaps $\g_j^+,$ $j=0,\ldots,q.$ We denote the set of bound states of $J$ by $\s_{\rm bs}(J).$

  Recall  that we consider Dirichlet boundary
condition $y_0=0$ in (\ref{pert}). Thus, any possible non-zero solution of
$Jy=\l y$ must have $y_1\neq 0,$ which implies that each eigenvalue of $J$ is simple (or else a linear combination would vanish
at $n=1$ and thus for all $n\in\N$).

 Moreover, if $f,g\in\ell^2_{\rm comp},$ where
$\ell_{\rm comp}^2$ denotes the $\ell^2$ functions on $\N$ with
a finite  support,  then the function    $\langle R(\l)f,g\rangle$
has an analytic extension from   $\L_1$ into the Riemann  surface
$\Lambda.$

\begin{definition}
\label{defstates} 1) A number $\l_0\in\L_2$ is a resonance if the
function $\langle R(\l)f,g\rangle$ has a pole at $\l_0$ for some
$f,g\in\ell^2_{\rm comp}.$ The set of resonances is denoted $\s_{\rm r}(J).$ The multiplicity of the resonance is
the multiplicity of the pole. If $\Re \l_0 =0,$ we call $\l_0$
un antibound state.\\
2) A real number $\l_0$ such that $\Delta^2(\l_0)=1$   is a virtual
state if  $\langle R(\l)f,g\rangle$ has a singularity at $\l_0$ for
some   $f,g\in\ell^2_{\rm comp}.$ The set of virtual
states is denoted $\s_{\rm vs}(J).$ \\
 3) The state  $\l_0\in\L$  is a bound state or a resonance or
 a virtual state of $J.$

 We denote the set of all states of $J$ by $\s_{\rm st}\,
 (J)=\s_{\rm bs}\,(J)\cup\s_{\rm r}\,(J)\cup\s_{\rm vs}\,(J)\ss \L.$

\end{definition}

The unperturbed Jacobi operator $J_0$ has one simple
state $\l_j$ in each   $\g^{\rm c}_j= \overline{\g^+_j}\cup
\overline{\g^-_j},$ $j=1,\ldots,q-1$ (see Proposition \ref{Prop_states_unperturbed}).
 Here the projection of $\l_j$ to $\C$ coincides with $\wt\l_j=\mu_j,$
 the zero of $\vp_q.$

Introduce the Jost solutions  $f^\pm=(f^\pm_n)_0^\iy$ and the
fundamental solutions $\vt^+=(\vt_n^+)_0^\iy, \vp^+=(\vp_n^+)_0^\iy$
to the equation
$$
{a}_{n-1} y_{n-1}+{a}_{n}y_{n+1}+{b}_ny_n=\l y_n,\qqq n\ge 1,
$$
 under the conditions
\[
\lb{defftf} f_n^\pm=\p_n^\pm,\qqq \vt_n^+=\vt_n,\qqq
\vp_n^+=\vp_n,\qqq n\geq p+1.
\]
Here and below $a_0=a_0^0=a_q^0$.  All functions $\vt_n^+,\vp_n^+,
n\ge0$   are polynomials. We rewrite the Jost solutions $f^\pm_n$ in
the form
\[
f^\pm_n=\vt_n^++m_\pm\vp_n^+,\qqq n\ge0.
\]
Note that for $\l\in\L_1$ we have $f^+(\l)\in\ell^2, $ and
$f^-(\l)=\ol{f^+(\ol\l)}.$  The functions
$f_n^\pm$ and the Titchmarch-Weyl functions $m_\pm$ are meromorphic
functions on $\L.$  Recall that the S-matrix for $J,J^0$ is   given
by
\begin{equation}
\lb{smatrix}
S(\l)={\ol{f_0^+(\l)}\/f_0^+(\l)}={{f_0^-(\l)}\/f_0^+(\l)} \qq {\rm
for} \qq \l \in \s_{\rm ac}(J^0).
\end{equation}

We pass to the formulation of main results of the  paper.
 Recall that if $\l\in\s_{\rm st}(J^0)$ then $\wt\l=\mu_j\in [\l_j^-,\l_j^+]$ for some $j=1,...,q-1,$ where $\mu_j$ denotes the Dirichlet eigenvalue and $\vp_q(\mu_j)=0$. Here the projection $\,\wt{}\,$ was introduced in (\ref{projection}). We describe all states of $J$.

\begin{theorem}\lb{th-charact-states}

\no i) The set of all state of $J$ has the decomposition
\[
\lb{Jdec}
\s_{st}(J)=\s^0(J)\cup \s^1(J),
\]
where
$$
\s^0(J)=\{\l\in \s_{st}(J^0): \vp_0^+(\wt\l)=0\}, \qqq \s^1(J)=\{\l\in \L: \l \not\in \s_{st}(J^0), f_0^+(\l)=0\}.
$$
 Moreover, each  $\l_0\in \s^0(J)$  is a simple state  of $J$
 and $0<|f_0^+(\l_0)|<\iy$.

\no ii) If $\l_1\in\L_1$ is a bound state of $J$, then
$\l_2\not\in\sigma_{\rm st}(J),$ where $\l_2\in\Lambda_2$ is the
same number as $\l_1$ but on the second sheet.

\no iii) Let $\l_0\in \L$  be a zero of $f_0^+$. Then $\vp_0^+(\wt\l_0)\ne 0$
 \end{theorem}

{\bf Remark.} 1) The proof of Theorem \ref{th-charact-states}
is given in Section \ref{s-genJost}.

\no 2) A state $\l_0\in \s^0(J)$  (bound, antibound or
virtual state)  is  not a zero of the Jost function $f_0^+$.
Moreover, $\l_0$ is a simple state of both $J$ and $J^0$. Such a
state is a singularity of the resolvent, but it is not a singularity
of the $S$--matrix (\ref{smatrix}).

In accordance with the continuous case \cite{KS} we define the
important function
\[\lb{FF}
 F(\l)=\vp_q(\l)f_0^+(\l)f_0^-(\l),\qqq \l\in \L_1.
\]
For the perturbation $V$ with $(u,v)\in\gX_\nu$ we define the
constants
\[
 \lb{c2_2p-1}
 c_3=c_1c_2,\qqq
 c_1={1\/\prod_{0}^p{a}_j},\qqq \qqq
 c_2=\ca c_1u_p(a^0_p+a_p) &\qq \mbox{if}\qqq \nu=2p,\\
c_1(a^0_p)^2{v}_p          & \qq\mbox{if}\qqq\nu=2p-1. \ac
\]
 The distribution of the states is summarized in
the following theorem.

\begin{theorem}
\label{mainstatictheorem}
Let  the  Jacobi operator $J=J^0+V$
satisfy \er{J0}--\er{V}. Suppose $(u,v)\in \gX_\n$, where $\n\in \{2p,2p-1\}$. Then the following facts hold true.\\
1) The function $F(\l),\l\in \L_1,$ is a real polynomial. Each zero
of $F$ is the projection of a state of $J$ on the first sheet. There
are no other zeros. Moreover, $F$ satisfies
\[
\lb{tn} F(\l)=-a_0^0\l^{\k}(c_3+\cO(\l^{-1})),\qqq  \k=\n+q-1,\qqq
\l\to \iy,
\]
here $\k$ is a total number of states (counted with
multiplicities).\\
2) The total number of bound and virtual states is  $\geq 2.$\\
3) In each finite open ``gap'' $\g_j^{\rm c}=\overline{\g_j^+}\cup
\overline{\g_j^-},$ $j=1,\ldots,q-1,$ there is always an odd
number $\ge 1$ of states (counted with multiplicities).\\
4) Let $\l_1<\l_2$ be any two bound states of $J,$ such that $\l_1,
\l_2\in\g_j^+,$ for some $j=0,\ldots, q.$ Assume that there are no
other eigenvalues on the interval
$\Omega^+=(\l_1,\l_2)\subset\g_j^+.$ Then there exists an odd number
$\geq 1$ of antibound states  on $\Omega^-,$ where
 $\Omega^-\ss\g_j^-\ss\L_2$ is the same interval
but on the second sheet, each antibound state being counted according
to its multiplicity.

\no 5)  $(-1)^{q-j}\dot{F}(\l)<0$ for any $\l\in\g_j^+\cap \s_{\rm
bs}(J), $  $j=0,1,\ldots,q$.\\
\no 6) If $\l\in\s_{\rm bs}(J)\cup\s_{\rm vs}(J)\cup\s^0(J)$,
 then $\l$ is a simple state of $J$.\\
\end{theorem}
The proof of Theorem \ref{mainstatictheorem}   follows from Lemmata \ref{Tm}--\ref{l-general case}.

{\bf Remark.} 1) The pre-image of a zero of $F$ is an eigenvalue or
a virtual state or a resonance of $J$. Thus we reformulate the problem
for the resolvent on the Riemann surface $\L$ as the problem for the
polynomial $F$ on the plane.

2)  There is an even number of
non-real resonances  since the resonances are zeros of the real polynomial $F$.

3) Due to this Theorem for the operator $J$ we define the vector-state
$\gr=(\gr_n)_1^\kappa$ by
\begin{align}
&\{\gr_j\}_{j=1}^\k=\s_{\rm st}(J),\qqq \gr_j\in \cup_0^q \g_n^+\in \L_1, \qqq  \gr_1<\gr_2<...<\gr_N, \qq
N\ge 0,\nonumber\\
&
\gr_j\in \L_2,\qqq 0\le
|\gr_{N+1}|\le|\gr_{N+2}|\le...\le|\gr_{\kappa}|,\lb{gr}
\end{align}
and the components of $\wt\gr$ are repeated according to the multiplicities of $\wt\gr_j $ as a zero of
the polynomial (\ref{FF}). Here $N$ is the number of bound states of $J$.

Now we pass to the inverse resonance problem.  We use the
parametrization $({u} ,{v} )=({u} _n,{v} _n)_1^p\in \R^{2p}$ for the
perturbation $V$ of the periodic coefficients of $J^0.$ We suppose that all  gaps are open:
$\lambda_j^-<\lambda_j^+,$ $j=1,\ldots,q-1$. We define the class of
all Jost functions on the Riemann surface $\L$ as follows.
\begin{definition}\lb{classJ}
For $\nu\in\N,$ let $\gJ_\nu$  denote the
class of rational functions  $f$  on $\L$ of the form
\begin{align*}
& f=P_1+m_+P_2, \\
& f(\l)=\ca c_1A_p+{\cO}(\l^{-1}) \qqq &{\rm if}
\qqq   \l\in\L_1\\
 -\frac{c_2 }{A_p}\l^{\nu}+\cO(\l^{\nu-1}) \qqq &{\rm
if}\qqq \l\in\L_2 \ac \qqq{\rm as}\qq \l\to \iy,
\end{align*}
where $c_1 >0,$ $c_2\ne 0$ and $P_1$ and $P_2$ are real polynomials (with real coefficients) of
the orders $\nu -2$ and $\nu -1$ respectively. Here
\[
\lb{defAp} A_p=\prod_{j=0}^pa_j^0.
\]
Let $\s(f)$ be the set of all zeros of $f$ on $\L$ and denote $\displaystyle \s_{\rm st}(f)=\s(f)\cup\s^0(f)\subset\L,$ where $$\displaystyle\s^0(f)=\{\l\in\s_{st}(J^0):\,\, P_2(\wt\l)=0\}.$$
We suppose that each zero of $f(\cdot)$ on the first sheet $\Lambda_1$ is real and belongs to
$\cup_0^q\ol\g_j^+.$ Let $$\displaystyle \s_{\rm
bs}(f)=\s_{\rm st}(f)\cap \cup_0^q\g_j^+.$$
Define the polynomial  $P(\l)
=\vp_q(\l)f(\l)f_-(\l),$ $\l\in \L_1$, where $f_-=P_1+m_-P_2$.\\
We suppose that the following properties hold true:\\
\no i)  if $\l\in\s(f)$, then $P_2(\l)\neq 0$, i.e., $\s(f)\cap \s^0(f) = \es$,\\
\no ii)  $(-1)^{q-j}\dot{P}(\wt\l)<0$ for any $\l\in\g_j^+\cap
\s_{\rm
bs}(f), $  $j=0,1,\ldots,q,$\\
\no iii) if $\l\in\s_{\rm bs}(f)\cup\s_{\rm vs}(f)\cup\s^0(f)$, where $\s_{\rm vs}(f)=\s(f)\cap(\cup_0^q\l_j^\pm),$
 then $\wt\l$ is a simple zero of $P.$ \\

\end{definition}

 Let $(u,v)\in \gX_\nu.$ Then from Theorems \ref{th-charact-states}, \ref{mainstatictheorem} and asymptotics in Section \ref{s-asym} it follows that
 the Jost function $f_0^+\in\gJ_\nu$ with $P_1=\vt_0^+,$  $P_2=\vp_0^+$ and $\s_{\rm st}(f_0^+)=\s_{\rm st}(J),$ $\s^0(f_0^+)=\s^0(J).$

Now we construct the mapping $\mF: \gX_\nu\to \gJ_\nu,$
$\nu\in\{2p-1, 2p\},$ by the rule:
\[
 ({u} ,{v} )\to f_0^+,
 \]
  i.e.  to each
$({u} ,{v} )\in \gX_\nu$ we associate $f_0^+\in \gJ_\nu$.

Our main inverse result is formulated in the following theorem.
\begin{theorem}\label{th-inverse}
 The mapping $\mF : \gX_\nu\to \gJ_\nu$ is one-to-one and onto. Moreover, the reconstruction algorithm is specified.
\end{theorem}

In Theorem \ref{th-inverse} we solve the inverse problem for mapping
$\mF.$ The solution is divided into the following three parts.
\begin{enumerate}
\item Uniqueness. Does the Jost function $f_0^+\in\gJ_\nu$ determine uniquely
 $({u} ,{v} )\in\gX_\nu$?
\item Reconstruction. Give an algorithm for recovering $({u} ,{v} )$ from
$f_0^+\in \gJ_\k$ only.
\item Characterization.  Give necessary and sufficient conditions for
$f_0^+$ to be the Jost functions for some  $({u} ,{v} )\in\gX_\nu.$
    \end{enumerate}

From Theorem \ref{th-inverse} it follows that  any $f\in\gJ_\nu$  is the Jost function $f_0^+$ for unique $J$ with $({u} ,{v} )\in\gX_\nu,$ and $P_1=\vt_0^+,$ $P_2=\vp_0^+,$ with the asymptotics
\begin{equation}\lb{as_vp_vt}
\vt_0^+=\frac{2a_0^2c_2}{A_p}\l^{\nu-2}+{\mathcal
O}\left(\l^{\nu-3}\right),\qqq
\vp_0^+=-\frac{2a_0c_2}{A_p}\l^{\nu-1}+{\mathcal O}
\left(\l^{\nu-2}\right),
\end{equation}
where  $c_2\neq 0,$   $A_p$ is given in (\ref{defAp}) and
$a_0=a_0^0=a_q^0$.

 Now we pass to the problem of reconstruction of the Jost function $f_0^+$ from  $\s_{\rm st}(J).$ Recall that $\s_{\rm st}(J)$  consists of the zeros of $f_0^+$ on $\L$ and the set $\sigma^0(J)$ (see Remark 2) after Theorem \ref{th-charact-states}).

 By Theorem \ref{mainstatictheorem}, 1),  the zeros of the polynomial $F$  defined in ( \ref{FF}) are given by $\{\wt\gr_j\}_{j=1}^\k,$ where
the set  $\{\gr_j\}_{j=1}^\k=\s_{\rm st}(J)$ satisfies (\ref{gr}). The polynomial $F$  can be uniquely reconstructed from the projection of all states $\{\wt\gr_j\}_{j=1}^\k$ and the constant $c_3$ in (\ref{tn}).

We have the following result.
\begin{theorem}\lb{th_getJost_1}
Suppose that $(u,v)\in \gX_\n$ and the polynomial $F$ has only simple zeros.  Then the Jost function $f_0^+$ is uniquely
determined by the polynomials
  $F$ and $\vp_0^+$.
 \end{theorem}

Now the polynomial $\vp_0^+$ can be reconstructed from its zeros and  the constant $c_2$ in (\ref{as_vp_vt}). Note that simple
examples show that zeros of the polynomial $\vp_0^+$ can be real and
non-real.

We have the identity
\begin{equation}
\lb{P_2}
 \vp_0^+=\frac{\vp_q}{2i\Omega}\left(f_0^+-f_0^-\right)
=\frac{\vp_q}{2i\Omega}f_0^+\left(1-S\right).
\end{equation}
Thus  the zeros of $\vp_0^+$ (under the conditions $\vp_q\neq 0$ and $\Omega\neq 0$) coincide with    the zeros of the function $1-S(\l)$ on $\L_1$  (see Lemma \ref{S=1}) and their multiplicities agree.

More precisely, let $Zeros\,(S-1)\in\L_1$  denote the set  of all zeros of
$S(\l)-1$ on $\L_1$ (counting the multiplicities).  Let  $\mu_j\in \overline{\gamma_j^+}\subset\L_1,$ $\vp_q(\mu_j)=0,$  $j=1,...,q-1,$ denote the
Dirichlet eigenvalue of $J_0.$

 From Lemma \ref{S=1} it follows that, if
 \begin{equation}\lb{hypS}
\left[Zeros\,(S-1)\setminus\left(\{\mu_j\}_{j=1}^{q-1}\cap\{\l_k^\pm\}_{k=0}^q\right)\right]\bigcap\left(\{\mu_j\}_{j=1}^{q-1}\cup\{\l_k^\pm\}_{k=0}^q\right)=\emptyset,
\end{equation}
 then the set $Zeros\,(S-1)$
 is the set of all zeros of $\vp_0^+.$
 We have the following result

 \begin{theorem}\lb{th_getJost_2}
 Suppose that the set of zeros $Zeros\,(S-1)$ on the first sheet $\L_1,$
   satisfy (\ref{hypS}), and  each zero of polynomial $F$ is simple.
Then the Jost function $f_0^+$ is uniquely determined by the polynomial $F,$  the set $Zeros\,(S-1)$ and  the constant  $c_2$.
\end{theorem}

Theorems \ref{th_getJost_1} and \ref{th_getJost_2} are proved in Section \ref{s-invres}.\\

{\bf Historical remarks.} A lot of papers is devoted to the
resonances for the Schr{\"o}dinger operator $-\frac{d^2}{dx^2}+q(x)$
on the line $\R$ and the half-line $\R_+$ with compactly supported
perturbation, see  \cite{Fr},  \cite{K4},  \cite{K5},
\cite{S},\cite{Z}, \cite{Z1},  and the references given there. Zworski [Z]
obtained the first results about the distribution of resonances for
the Schr\"odinger operator with compactly supported potentials on
the real line. One of the present authors obtained the uniqueness,
the recovery and the characterization of the $S$-matrix for the
Schr\"odinger operator with a compactly supported potential on the
real line \cite{K4} and on the half-line \cite{K5}, see also
\cite{Z1}, \cite{BKW} concerning the uniqueness.

The problem of resonances for the Schr{\"o}dinger with periodic plus
compactly supported potential $-\frac{d^2}{dx^2}+p(x)+q(x)$  is much
less studied: \cite{F1}, \cite{KM}, \cite{K1}, \cite{KS}. The
following results were obtained in \cite{K1}, \cite{KS}: 1) the
distribution of resonances in the disk with large radius is
determined, 2) some inverse resonance problem, 3) the existence of a
logarithmic resonance-free region near the real axis.
Note that in our paper we use the methods from  \cite{KS}, modified for the Jacobi operator $J$.

Finite-difference Schr{\"o}dinger and Jacobi operators express many
similar features. Spectral and scattering properties of infinite
Jacobi matrices are much studied (see \cite{Mo}, \cite{DS1},
\cite{DS2} and references given there). The inverse problem  for
periodic Jacobi operators $J^0$ was solved in \cite{BGGK},
\cite{K3}, \cite{KKu}, \cite{Mo}, \cite{P} and see
 references therein.

The inverse resonances problem was recently solved in the case of
constant background \cite{K2}.
 The inverse scattering problem  for asymptotically periodic
 coefficients was solved by  Khanmamedov: \cite{Kh1}
 (on the line, note that the russian versions were dated much earlier), \cite{Kh2}
 (on the half-line) and Egorova, Michor, Teschl \cite{EMT}
 (on the line in case of quasi-periodic background).


In our paper we apply some results from \cite{Kh1}, \cite{Kh2} and
\cite{EMT}. There were some mistakes in the paper \cite{EMT},
\cite{BE}. Some of them we correct in   Section
\ref{s-periodic}. However, in our context  of finite rank
perturbations their results still hold in the original form.

We plan to apply the results of  our paper to the Schr{\"o}dinger operator
on nanotubes (see \cite{IK1} and references therein). The similar methods
are applied in \cite{IK2} and \cite{IK3} to the direct and the inverse resonance
problems on the line.

{\bf Plan of the paper.}
 In Part \ref{p-direct} we consider the
direct problems for the Jacobi operators on the half-line. In
Section \ref{s-periodic} we recall some well known facts about the
periodic Jacobi operators and describe the states for the periodic
Jacobi operators on the half-line. We present also the revised construction of the
quasi-momentum map. In Section \ref{s-genJost} we
consider the properties of the  Jost functions and prove  Theorems
\ref{th-charact-states} and \ref{mainstatictheorem}.

Part \ref{s-inverse} is devoted to the inverse resonance problem.
 In Section \ref{s-scatt} we recall the results
of Khanmamedov on the inverse scattering problem on the half-line
which we apply in Section \ref{s-invres} and prove the inverse
 results.

In Part \ref{s-asym} we collect the asymptotics of the Jost functions
which we need in the proofs.

{\bf Acknoledgement.} The authors are  indebted to the referee for numerous comments and suggestions.

\section{Direct problem}\lb{p-direct}
\setcounter{equation}{0}

\subsection{Unperturbed Jacobi operators $J^0$. }
\lb{s-periodic} We need some known properties of the $q-$periodic
Jacobi operator $J^0$ on $\N$ (see \cite{P}, \cite{T}, \cite{Kh1}).
Recall that the fundamental solutions $\vp=(\vp_n)_0^\iy$ and $
\vt=(\vt_n)_0^\iy$ and the Lyapunov function $ \D$ were defined in the
Introduction.  The spectrum of $J^0$ consists of an absolutely
continuous part $\s_{\rm ac}(J^0)= \bigcup\limits_{1}^q \s_j$ plus
at most one eigenvalue in each non-empty gap $\g_j,$ $ j=1,...,q-1$,
 where the bands $\s_j$ and the gaps $\g_j$ are given by
\er{ends}.

If there are exactly $N\ge 1$ nondegenerate gaps in the spectrum of
$\s_{\rm ac}(J^0)$, then the operator $J^0$ has exactly $N$ states;
the closed gaps $\g_n=\es$ do not contribute to any states. In
particular, if all $\g_j=\es, j\ge 1$, then $q=1$  (see \cite{BGGK},
\cite{KKu}, \cite{K3}) and $J^0$ has no states.  A more detailed
description of the states of $J^0$ is given in Proposition
\ref{Prop_states_unperturbed} below.

 In each finite ``gap'' $[\l_j^-,\l_j^+],$ $j=1,\ldots,q-1,$ there is
 one simple zero of
polynomials $\vp_q(\l),$ $\dot\D(\l),$ $\vt_{q+1}(\l)$. Here $\l_1^\pm,\ldots,\l_{q-1}^\pm$ are all
endpoints of the bands, see \er{ends}. Note that
$\Delta (\l_j^\pm)=(-1)^{q-j}$. The sequence of zeros
 of the polynomial $\Delta^2-1$ of degree $2q$
can be enumerated by $(\l_j^\pm)_0^q,$ $\l_0^+=\l_0^-,$ $\l_q^+=\l_q^-.$ We have
\begin{align*}
&\vp_q={a_0^0}\prod_{j=1}^{q-1}(\l-\mu_j),\qqq
\vt_{q+1}=-a_0^0\prod_{j=1}^{q-1}(\l-\nu_j),\\
& \D^2-1=\frac{1}{4}(\l-\l_0^+)(\l-\l_q^-)\prod_{j=1}^{q-1}(\l-\l_j^-)(\l-\l_j^+),
\end{align*}
where $\mu_j\in [\l_j^-,\l_j^+]$ are  the zeros of ${\vp_q}$ and
$\nu_j\in [\l_j^-,\l_j^+]$ are  the zeros of ${\vt_{q+1}}$
(Dirichlet or Neumann eigenvalues).  We put
 $$
 A=A_q=\prod_{j=1}^qa_j^0=1,\,\qqq B=\sum_{j=1}^qb_j^0.
 $$
Note the following asymptotics:
\[\lb{as_z}
\vp_q=a_0^0\l^{q-1}+\cO(\l^{q-2}),\,\qqq
\D(\l)={z^q+z^{-q}\/2}={\l^q\/2}+{\mathcal
O}(\l^{q-1})\,\,\mbox{as}\,\,\l\rightarrow\iy.
\]
Here the function $z=z(\l)=e^{i\vk(\l)}$ is explained later  in this section and  $\vk(\l)$ is the quasi-momentum satisfying (\ref{ca}).

Recall that $\G$ is the complex $\l$-plane with cuts along the segments $\s_j,$
$j=1,2,\ldots,q.$
$\G$ will be identified with the first sheet $\Lambda_1$.
 We use the standard definition of the root:
$\sqrt{1}=1$ and fix the branch of the function $\sqrt{\D^2(\l)-1}$
on $\L$ by demanding $\sqrt{\D^2(\l)-1} <0$ for $\l>\l_q^-, \l\in
\G$ (in accordance with (\ref{branch})). We define the  first $\x_+$ and
the second $\x_-$ Floquet multipliers on the plane $\L_1$ or $\G$  by
$$
\x_\pm(\l)=\Delta (\l) \pm\sqrt{\Delta^2(\l)-1},\qqq \l\in \L_1.
$$
By our choice of the branch we have $|\x_+(\l)|<1,$  $|\x_-(\l)|>1$ and
\[
\lb{Dfactor}
\sqrt{\Delta^2(\l)-1}=-\frac{1}{2}\sqrt{\l-\l_0^+}\sqrt{\l-\l_q^-}\,\prod_{j=1}^{q-1}\sqrt{\l-\l_j^-}\sqrt{\l-\l_j^+}.
\]
for all $\l\in \L_1$.  The functions $\x_\pm(\l)$ are continuous up to the boundary
$\partial\L_1$  and
$|\x_\pm(\l)|=1\,\,\mbox{for}\,\,\l\in\partial\L_1.$ Moreover for
$\l\in\Lambda_1,$
$$
\xi^\pm(\l)=(2\Delta(\l))^{\mp 1}\left(1+{\mathcal
O}\left(\l^{-2q}\right)\right)=\l^{\mp q}\left( 1\pm
\frac{B}{\l}+\cO\left(\frac{1}{\l^2}\right)\right).
$$

For two sequences $x=(x_n)_1^\iy, y=(y_n)_1^\iy$  we introduce the
unperturbed  Wronskian by
\[
\lb{W} \{x,y\}_n^0={a}_n^0(x_ny_{n+1}-x_{n+1}y_n).
\]

Using that $\{x,y\}_n^0$ is independent of $n$ for two solutions of (\ref{1e2}) and putting $x=\vt,$ $y=\vp,$ we apply the conditions $\vt_0=\vp_1=1,$ $\vt_1=\vp_0=0$ and obtain
\[
\lb{Fact7} 1-\D^2+\f^2=1-\vp_{q+1}\vt_{q}=-\vp_{q}\vt_{q+1}.
\]
Thus, we get
\[
\lb{mm}
m_+m_-=-\frac{\vt_{q+1}}{\vp_q}.
\]
 This identity considered at zeros of polynomial $\vp_q(\l)$ of degree $q-1$
 shows: if one of the solutions $\psi_n^\pm(\l)$ is regular, then the other has
simple poles, one in each
 finite gap $\gamma_n,$ $n=1,\ldots,q-1$.

 Equation (\ref{J0}) has two Bloch solutions
$\psi_n^\pm=\psi_n^\pm(\l)$ which satisfy $\psi_{kq}^\pm=\xi_\pm^k,$
$k\in\Z,$ and at the end points of the gaps we have
$|\psi^\pm_{kq}(\l_n^\pm)|=1.$ As for any $\l\in\L_1$ we have
$\psi^+\in \ell^2(\N),$ then functions $\psi^\pm(\l)$ are the
Floquet solutions for (\ref{J0}).

Now we consider the spectrum of the half-infinite Jacobi matrix
$J^0$ defined by (\ref{Jm}) or (\ref{pert}) with coefficients
${a}^0_j,$ ${b}^0_j,$ $j\in\N,$ verifying (\ref{pert}).

\begin{proposition}[States of $J^0$]\lb{Prop_states_unperturbed}
The unperturbed operator $J^0$
 has absolutely continuous spectrum (\ref{acspectr}):
$\s_{\rm ac}(J^0)=\cup_{j=1}^q\s_j$ and one simple
state $\l_j$ in each   $\g^{\rm c}_j= \overline{\g^+_j}\cup
\overline{\g^-_j},$ $j=1,\ldots,q-1.$    Here the projection of $\l_j$ on $\C$ coincides with $\wt\l_j=\mu_j,$
 the zero of $\vp_q.$
\end{proposition}
{\bf Proof.}
The kernel of the resolvent of $J^0$  is given by
$$
R^0(n,m)=-\frac{\varphi_n\psi_m^+}{\left\{\varphi,\psi^+\right\}}=
\frac{\varphi_n\psi_m^+}{a_0^0},\,\,n <m,
$$
since $\left\{\varphi,\psi^+\right\}=-a_0^0.$ According to Lemma
\ref{defstates2} (see Section \ref{s-genJost}), the bound states
(resonances) are the poles of
${\mR}^0_n=\psi^+_n(\lambda)=\vt_n(\lambda)+m_+(\l)\varphi_n(\l)$ or
of $m_+(\l)$
 on $\Lambda_1$ (respectively on $\Lambda_2$).

 From
(\ref{mm}) it follows that if  $\mu_n\neq\l_n^\pm,$  $n=1,\ldots,q-1,$  then
 one of the following two cases holds true:

 (i) $m_+$ has simple pole at
$\mu_n,$  $m_-$ is regular and $\mu_n$ is the bound state,

(ii) $m_-$ has simple pole at $\mu_n,$ $m_+$ is regular and
$\mu_n$ is the antibound state.

Now suppose that either $\mu_n=\l_n^-,$ $\l_0=\mu_n +\epsilon$
or $\mu_n=\l_n^+,$ $\l_0=\mu_n -\epsilon,$ $\epsilon >0.$
   Then \[\lb{unpertvirt} m_+(\l_0)=\frac{c}{\sqrt{\epsilon}}
   +{\mathcal
O}(1),\,\,\epsilon\rightarrow 0,\,\,c\neq 0.
\]
 Moreover, for $n\neq 0,q,$
$\psi^+_n(\l_0)=\vt_n(\mu_n)
+\left(\frac{c}{\sqrt{\epsilon}}+{\mathcal
O}(1)\right)\varphi_n(\mu_n),$  the function $(\mR_n^0(.))^2$
 has a pole at $\mu_n$ for almost all $n\in\N$ and $\mu_n$ is the virtual state.
 \hfill\BBox

We have also
$$
m_\pm=\frac{\x_\pm-\vt_q}{\vp_q}.
$$
Moreover, $\mu_j\in\g_j$ is the antibound state iff $ \x_+(\mu_j)=\vt_q
(\mu_j)$ and $\mu_j\in\g_j$ is the bound state iff $\x_-(\mu_j)=\vt_q (\mu_j).$ Note that
on each $\g_j^+,$ $j=0,1,\ldots,q,$ $m_\pm$ are real functions.

{\bf Quasi-momentum map and Riemann surface $\cZ.$}\lb{ss-cut}

 We construct the conformal mapping  of the Riemann surface onto
 the plan with ``radial cuts'' $\cZ.$   Our definition corrects
 the similar construction  in \cite{BE} and \cite{EMT}, where there
 was a mistake.

We suppose that all  gaps are open:
$\lambda_j^-<\lambda_j^+,$ $j=1,\ldots,q-1$.

Introduce a domain $\C\sm\cup_0^{q}\ol\g_j$ and a quasi-momentum domain
$\K$ by
$$
\K=\{\vk\in \C:-\pi\le\Re \vk\le 0\}\sm\cup_1^{q-1}\ol\G_j,\ \ \G_j =\rt(
-{\pi j+ih_j\/q}, -{\pi j-ih_j\/q}\rt).
$$
Here $h_j\ge 0$ is defined by the equation $\cosh
h_j=(-1)^{j-q}\D(\alpha_j)$ and $\alpha_j$ is a zero of $\D'(\l)$ in
the ``gap'' $[\l_j^-,\l_j^+]$. For each periodic Jacobi operator
there exists a unique conformal mapping $\vk:\C\sm\cup_0^{q}\ol\g_j\to
\K$ such that the following identities and asymptotics hold true:
\[
\lb{ca}
  \cos q\vk(\l)=\D (\l),   \ \ \ \l\in \C\sm\cup_0^{q}\ol\g_j,
  \ \ \ \ \ {\rm and} \ \ \
  \ \ \ \vk(it)\to \pm i \iy\ \ \ {\rm as} \ \ t\to \pm\iy.
\]
The quasi-momentum $\vk$ maps the half plane
$\C_\pm=\{\l\in\C;\,\,\pm\Im\l>0\}$ onto the half-strip
$\K_\pm=\K\cap \C_\pm$ and $\s_{\rm
ac}(J^0)=\{\l\in\R;\,\,\Im\vk(\l)=0\}$.

Define the two strips $\K_S$ and $\cK$ by
$$
\K_S=-\K\qqq \mbox{and} \qqq \cK=\K_S\cup \K\ss
\{\vk\in \C: \Re \vk\in [-\pi,\pi]\}.
$$
The function $\vk$ has an analytic continuation from $\L_1\cap \C_+$
into $\L_1\cap \C_-$  through the infinite gaps $\g_q=(\l_q^-,\iy)$
by the symmetry and satisfies:

1)  $\vk$  is a  conformal mapping $\vk:\L_1\to \cK_+=\cK\cap \C_+$,
where we identify the boundaries $\{\vk=\pi+it, t>0\}$ and
$\{\vk=-\pi+it, t>0\}$.

2) $\vk:\L_2\to \cK_-=\cK\cap \C_-$  is a  conformal mapping, where
we identify the boundaries $\{\vk=\pi-it, t>0\}$ and $\{\vk=-\pi-it,
t>0\}$.

3) Thus $\vk:\L\to \cK$  is a  conformal mapping.

Consider the function $z=e^{i\vk(\l)},\,\, \l \in \L$. The function $z(\l),$ $\l \in \L,$
is a conformal mapping $z:\L\to \cZ=\C\sm \cup \ol g_j$, where the radial cut $g_j$ is given by
$$
g_j=(e^{-{h_j\/q}+i{\pi j\/q}}, e^{{h_j\/q}+i{\pi j\/q}}), \qqq j=\pm 1, ..., \pm (q-1).
$$

The function $z(\l),$  $\l \in \L,$ maps the first sheet $\L_1$ into the ``disk'' $\cZ_1=\cZ\cap \dD_1,$
$\dD_1=\{z\in\C:\,\,|z|<1\},$
and  $z(\cdot)$ maps the second sheet $\L_2$ into the domain  $\cZ_2=\cZ\sm \dD_1$.
In fact, we obtain the parametrization of the two-sheeted Riemann surface $\L$ by the ``plane''
$\cZ$. Thus below we call $\cZ_1$ also the ``physical sheet'' and $\cZ_2$ also the ``non-physical sheet''.

Note that if all $a_n^0=1, b_n^0=0$, then we have $\l={1\/2}(z+{1\/z})$. This function
$\l(z)$ is a conformal mapping from the disk $\dD_1$ onto the cut domain $\C\sm [-2,2]$.

Now, the functions $\psi^\pm(\l)$ can be considered as functions of $z\in \cZ$.
The functions $\psi^\pm_n(z)\equiv \psi^\pm_n(\l(z))$ are meromorphic in $\cZ$ with the only possible singularities at the images of the Dirichlet eigenvalues $z(\mu_j)\in\cZ$ and at $0.$ More precisely,\\
\no 1) $\psi_n^\pm$ are analytic in $\cZ \sm (\{z(\mu_j)\}_{j=1}^{q-1}\cup\{0\})$ and continuous up to $\pa \cZ\sm\{z(\mu_j)\}_{j=1}^{q-1}.$

\no 2) $\psi_n^\pm(z)$ has a simple pole at $z (\mu_j)\in\cZ$  if $\mu_j$ is a pole of $m_\pm,$ no pole if $\mu_j$ is not a singularity of $m_\pm$ (not a square root singularity if $\mu_j$ coincides with the band edge)  and if $\mu_j$  coincides with the band edge:  $\mu_j=\l_j^\sigma,$ $\sigma=+$ or $\sigma=-$, $j=1,\ldots, q-1,$ then
\[
\lb{B2psi} \psi_n^\pm(z)=\pm  \sigma(-1)^{q-j}\frac{iC(n)}{z -z(\l_j^\sigma)} +\cO(1),\qq\l\in [\l_{j-1}^+,\l_j^-],
\]
for some constant  $C(n)\in \R$. Note that the sign comes from the analytic continuation of the square root $\Omega(\l)$ using the definition
(\ref{branch}).

\no \no 3) The following identities hold true:
\[\lb{prop4}
\psi_n^\pm(\overline{z})=\psi_n^\pm(z^{-1})= \psi_n^\mp(z)=\overline{\psi_n^\pm(z)}\,\,\mbox{as}\,\, |z| =1.
\]
\no 4) The following asymptotics hold true:
$$
\psi_n^\pm(z)= (-1)^n\lt(\prod_{j=0}^{n-1}{}^*a_j\rt)^{\pm 1}z^{\pm
n}\lt(1 +\cO (z)\rt)
 \qqq \mbox{as} \qqq z\rightarrow 0.
 $$

We collect below some properties of the quasi-momentum $\vk$ on the gaps.

 On each $\g_j^+, j=0,1,\ldots,q,$ the
quasi-momentum $\vk(\l)$  has  constant real part and positive
$\Im\vk$:
 $$
 \Re \vk|_{\g_j^+}= -\frac{q-j}{q}\pi,
 \qqq\vk(\l_j^-)=\vk(\l_j^+)=-\frac{q-j}{q}\pi,\qqq
 \Im\vk|_{\g_j^+}>0.
 $$
  Moreover, as $\l$ increases from $\l_j^-$ to $\a_j$
  the  imaginary part  $\Im\vk\equiv h(\l)$
is monotonically increasing from $0$ to $h_j$ and as $\l$ increases
from $\a_j$ to $\l_j^-$   the  imaginary part  $\Im\vk\equiv h(\l+i0)$
is  monotonically  decreasing from $h_j$ to $0$. Then
 \[
 \lb{isin}
\frac12\vp_q(\l)(m_+(\l)-m_-(\l))= \sqrt{\D^2(\l)-1}=i\sin q\vk(\l)= -(-1)^{q-j}\sinh qh(\l+i0),
\]
where $\sinh qh=-2^{-1}(z^q-z^{-q})>0.$

\

\subsection{ The perturbed Jacobi operator,  Jost functions.}\lb{s-genJost} We
consider the operator  $J=J^0+V$ given by (\ref{pert}). Recall that
$f_n^\pm$ are solutions to the equation
\[
\lb{eq-pert}
{a}_{n-1} y_{n-1}+{a}_{n}y_{n+1}+{b}_ny_n=\l y_n,\qq
  \l\in\L_1,
  \]
 satisfying
\[
\label{defJost} f_n^\pm=\p_n^\pm,\qqq \mbox{for all }\,\, n\geq p+1.
\]
Recall that ${a}_n={a}^0_n+{u}_n,\qq
 {b}_n= b_n^0+{v}_n$.
Equation (\ref{eq-pert}) has unique solutions $\vt_n^+,$ $\vp_n^+$
such that
$$
\vt_n^+=\vt_n,\qqq \vp_n^+=\vp_n,\qqq \mbox{for all }\,\, n\geq p+1.
$$
The functions $\vt_n^+(\cdot),$ $\vp_n^+(\cdot)$ are polynomials.
The functions $f_n^\pm$ have the form \begin{equation}\lb{fform}
f^\pm_n=\vt_n^+ +m^\pm
\vp_n^+\end{equation} and satisfy $\ol f^\pm_n(\ol{\l})=f^\mp_n(\lambda),$
$\l\in\G.$

\begin{lemma}
\lb{tph+} The zeros of the polynomials $\vt_0^+$ and $\vp_0^+$ are
disjoint.
\end{lemma}
{\bf Proof.} Assume that $\vt_0^+(\l_0)=\vp_0^+(\l_0)=0$ for some
$\l_0\in \C$. Then $\vt_n^+(\l_0)=a\vp_n^+(\l_0)$ for all $n\ge 1$
and some $a\ne 0$. Then \er{defftf} gives $\vt_n(\l_0)=a\vp_n(\l_0)$
for all $n> p$ and thus $\vt_n(\l_0)=a\vp_n(\l_0)$ for all $n> 1$
and the Wronskian $\{\vt(\l_0),\vp(\l_0)\}=0$. We have a
contradiction, since $\{\vt(\l_0),\vp(\l_0)\}=1$. \hfill\BBox

By Definition \ref{defstates} a state is a singularity of the resolvent.
 The
kernel of the resolvent of $J$ is given by
$$
R(m,n)=\langle e_m,(J-\l)^{-1} e_n\rangle= -\frac{\Phi_m
f_n^+}{\left\{\Phi,f^+\right\}}={\Phi_m {\mR}_n(\l)\/a_0},\,\,m <n,
$$
$$
{\mR}_n(\l)=\frac{f^+_n(\lambda)}{f^+_0(\lambda)}.
$$
Here $e_n=(\delta_{n,j})_1^\iy$ is the unit vector in $\ell^2$,  and $\F=(\F_n)_0^\iy$ is a
solution of the equation \er{eq-pert} under the condition $\Phi_0=0,$ $\Phi_1=1,$ and note that
$\left\{\Phi,f^+\right\}=-{a}_0f^+_0.$ Each function
$\Phi_n(\lambda),$ is polynomial in $\lambda.$ The
function $R(n,m)$ is meromorphic on $\Lambda$  for each $n,m\in\N.$
Then the singularities of $R(n,m)$ are given by the singularities of
${\mR}_n(\l)$. We have
\begin{lemma}\label{defstates2}
 1) A real number $\l_0\in\g_k^+,$ $k=0,1,\ldots,q$ is a bound state,
 if the
 function ${\mR}_n(\l)$ has a pole at $\l_0$ for some $n\in\N$. Recall (see Introduction, before Definition \ref{defstates}) that
 the bound states are simple.\\
 2) A number $\l_0\in\Lambda_2,$ is a resonance, if the
 function ${\mR}_n(\l)$ has a pole at $\l_0$ for some $n\in\N$. The
 multiplicity of the resonance is the multiplicity of the pole.\\
3) A real number $\l_0=\l_k^\pm,$ $k=0,\ldots,q,$  is a virtual
state if $\mR_n^2(\l)$ or $\mR_n(\l)$ has a pole at $\l_0$ for some
 $n\in\Z_+$.
\end{lemma}

{\bf Proof of Theorem \ref{th-charact-states}}
i) We start with the case $\l_0\not\in\sigma_{\rm st}(J^0).$

Let $\Omega(\l_0)\neq 0.$ Then $f_n^+,$ $n\in\N,$  is analytic at $\l_0\in\L.$ Then ${\mR}_n(\l)$ has a pole at $\l_0$ iff $f_0^+(\l_0)=0.$

Let now  $\Omega(\l_0)= 0.$  Using (\ref{Dfactor}) we get $
m^\pm(\l)= m^\pm (\l_0)+c\sqrt{\epsilon}+{\mathcal
O}(\epsilon),$ $\l-\l_0=\epsilon\rightarrow 0,$
and $c \neq 0.$ We distinguish between two cases.\\
a) Firstly, let $\varphi^+_0(\wt\l_0)\neq 0.$ Then identity
$f_0^+(\l_0)=\vt^+_0(\l_0) +m^+\varphi^+_0(\l_0)=0$  implies
(\ref{3.6})
$$f_0^+(\l)=\varphi^+_0(\wt\l_0)c\sqrt{\epsilon} +{\mathcal
O}(\epsilon),\,\,\mR_n(\l)=\frac{f_n^+(\l)}{\varphi^+_0(\wt\l_0)c
\sqrt{\epsilon}}(1+{\mathcal
O}(\sqrt{\epsilon})),\,\,c\varphi^+_0(\wt\l_0)\neq 0.$$ Then $\l_0$ is a virtual state of $J.$\\
b) Secondly, if $\varphi^+_0(\wt\l_0)=0,$ then we obtain
$\vt^+_0(\wt\l_0)\neq 0$ by  Lemma \ref{tph+} and
$f_0^+(\l_0)=\vt^+_0(\wt\l_0)\neq 0.$ Then $\l_0$ is not a singularity of the resolvent.

Now we consider the case $\l_0\in\sigma_{\rm st}(J^0).$  Then $\vp_q(\wt\l_0)=0.$\\
  Suppose firstly that $\Omega(\l_0)\neq 0.$ Then $\l_0$ is a pole of $m_+$ and therefore $\l_0$ is a pole of the Jost solution $f_n^+(\l)=\vt_n^++m_+\vp_n^+$ on either $\L_1$ or $\L_2$ for all $n\in\N$ such that $\vp_n^+(\wt\l_0)\neq 0.$
 Then using Lemma \ref{tph+} we get that if $\vp^+_0(\wt\l_0)=0$ then $f_0^+(\l_0)\neq 0$  and
 $\l_0$ is a pole of
$${\mR}_n(\l)=\frac{f_n^+(\l)}{f_0^+(\l)}=
\frac{\vt^+_n+m_+\vp^+_n}{\vt^+_0+m_+\vp^+_0},\qq n\in N,$$
iff $\vp^+_0(\wt\l_0)=0.$  Moreover, $\l_0$ is a simple state (as a pole of $m^+$).

Suppose now that $\l_0\in\sigma_{\rm st}(J^0)$ and $\Omega(\l_0)= 0.$\\
  Then we have (\ref{unpertvirt}):
$$m^+(\l)=\frac{c}{\sqrt{\epsilon}}+{\mathcal
O}(1),\,\,\l-\l_0=\epsilon\rightarrow 0,\,\,c\neq 0.$$ We distinguish between two cases.\\
a) Firstly, let $\varphi^+_0(\wt\l_0)\neq 0.$ Then identity
$f_0^+(\l_0)=\vt^+_0(\l_0) +m^+(\l_0)\varphi^+_0(\l_0)=0$ implies
$$
f_0^+(\l)=\frac{\varphi^+_0(\wt\l_0)
c}{\sqrt{\epsilon}}+\cO(1),\,\qqq
\frac{f_n^+(\l)}{f_0^+(\l)}=\frac{\vt_n^+(\wt\l)
+\left(\frac{c}{\sqrt{\epsilon}}+{\mathcal
O}(1)\right)\varphi^+_n(\wt\l)}{\frac{\varphi^+_0(\wt\l_0)
c}{\sqrt{\epsilon}}+\cO(1)}=\frac{1+{\mathcal
O}(\sqrt{\epsilon})}{\varphi^+_0(\wt\l_0)},
$$
and each function $\mR_n(.),$ $n\in\N,$ does not have singularity at
$\l_0.$\\
b) Secondly, let $\vp_0^+(\wt\l_0)=0.$
Then
$f_0^+(\l_0)=\vt^+_0(\wt\l_0)\neq 0$ by Lemma \ref{tph+}. Moreover, we obtain
$f^+_n(\l)=\vt_n^+(\wt\l) +\left(\frac{c}{\sqrt{\epsilon}}+{\mathcal
O}(1)\right)\varphi^+_n(\wt\l),$ and the function $(\mR_n(.))^2,$
$n\in\N,$ has  simple pole at $\l_0.$

ii) Suppose $\l_1\in\L_1$ is a bound state of $J$ and $\l_1\not\in\sigma_{\rm st}(J^0).$ Then by i) we have $f_0^+(\l_1)=0$ and as $\{f^+,f^-\}\neq 0$ we have $f_0^-(\l_1)\neq 0$ (by the argument similar to Lemma \ref{tph+}). The last identity is equivalent to $f_0^+(\l_2)\neq 0$ for  $\l_2\in\L_2$ such that $\wt\l_2=\wt\l_1.$

iii) In i) it was shown that if $\l_0\in \sigma_{\rm st}(J^0)$ then $f_0^+(\l_0)\neq 0.$  So it is enough to consider the case $\l_0\in\L$ is a zero of $f_0^+$ and $\l_0\not\in \sigma_{\rm st}(J^0).$  If $\vp_0^+(\wt\l_0)=0$ then $f_0^+(\l_0)=\vt_0^+(\wt\l_0)\neq 0$ as in ii) which is a contradiction.

\hfill\BBox

Define the function
\[
\lb{dfn} F_n(\l)=\vp_q(\l)f_n^+(\l)f_n^-(\l),\qq \l\in \L_1.
\]
Note that $F_0=F$ defined previously  in (\ref{FF}).
Using (\ref{fform}) and (\ref{TitchWeyl}), (\ref{Fact7}),  (\ref{mm}) we get
\[
\lb{aFn} F_n=\vp_q(\vt_n^+)^2+2\f
\vt_n^+\vp_n^+-\vt_{q+1}(\vp_n^+)^2, \qq n\ge 0.
\]
The following Lemma is proven in Section \ref{s-asym}.

\begin{lemma}
\lb{Tm} Let $\nu\in \{ 2p, 2p-1\}.$ Each function $F_n(\l)=\vp_q(\l)f_n^+(\l)f_n^-(\l),$ $n\ge
0,$ is a polynomial and satisfy
\[
\lb{aFn1}
 F_n(\l)=-a_0^0\l^{\k-2n}\lt(c_3(n)+\cO(\l^{-1})\rt),
 \qqq \k=\nu+q-1,\qqq \l\to \iy,
\]
\[
 \lb{c2_2p-1}
 c_3(n)=c_1(n)c_2(n),\qq
 c_1(n)={1\/\prod_{j=n}^p{a}_j},\qq
 c_2(n)=\ca c_1(n)u_p(a^0_p+a_p) &\qq \mbox{if}\qqq \nu=2p,\\
c_1(a^0_p)^2{v}_p          & \qq\mbox{if}\qq\nu=2p-1. \ac
\]
\end{lemma}
{\bf Remark.} It follows that the function
$F_n(\l)=\vp_q(\l)f_n^+(\l)f_n^-(\l)$ is polynomial of degree
$2(p-n)+q-1$ (if ${u} _p\neq 0$) or $2(p-n)+q-2$ (if ${u} _p=0,$
${v} _p\neq 0$).
 From  the asymptotics (\ref{F_n_1}),  (\ref{F_n_2}) collected in Section
\ref{s-asym}, we get the sign of $F=F_0$ as $\lambda\rightarrow\infty:$
$$
\sign F(\l)=\ca  -\sign u_p & if \ u_p\neq 0\\
                 -\sign({v}_p) & if \ a_p^0\neq{a}_p\ac \qq  as \ \l\to
                 \infty,
$$
$$
\sign F(\l)=\ca (-1)^{2p+q-2}\sign u_p &if\
u_p^0\neq 0\\
-(-1)^{2p+q-2} \sign({v} _p) &if \ u_p^0=0, v_p\neq 0\ac \qq as \
\l\to
                 -\infty.
$$

We summarize the results about the virtual states $\sigma_{\rm vs}(J)$ obtained in the proof of Theorem \ref{th-charact-states} in the following Lemma.

\begin{lemma}[Virtual states]\lb{L-virt}
Let $\l_0=\l_k^\pm$ for some $k=0,\ldots,q-1.$  If $\l_0=\l_k^+$
then put $\lambda=\l_0-\epsilon.$ If  $\l_0=\l_k^-,$ then put
$\lambda=\l_0+\epsilon.$ Here
$\epsilon>0$ is small enough.\\
i) Let $\l_0\not\in\sigma_{\rm st}(J^0)$  and  $f_0^+(\l_0)=0$.
Then $\wt\l_0$ is a simple zero of $F,$   $\l_0$ is virtual state of
$J$ and
\[\lb{3.6}
f_0^+(\l)=\varphi^+_0(\wt\l_0)c\sqrt{\epsilon} +{\mathcal
O}(\epsilon),\,\,\mR_n(\l)=\frac{f_n^+(\l)}{\varphi^+_0(\wt\l_0)c
\sqrt{\epsilon}}(1+{\mathcal
O}(\sqrt{\epsilon})),\,\,c\varphi^+_0(\wt\l_0)\neq 0.
\]
ii) Let  $\l_0\in\sigma_{\rm st}(J^0)$   and $\varphi^+_0(\wt\l_0)\neq
0$. Then $F(\wt\l_0)\neq 0$ and each $\mR_n(.),$ $n\in\N,$ does not
have singularity at $\l_0$ and $\l_0$ is
not a virtual state of $J.$\\
iii) Let $\l_0\in\sigma_{\rm st}(J^0)$  and $\varphi^+_0(\wt\l_0)=
0$. Then $\l_0$ is virtual state of $J,$ $f_0^\pm (\l_0)\neq 0,$
$\wt\l_0$ is simple zero of $F,$ and each $(\mR_n(.))^2,$ $n\in\N,$ has
pole at $\l_0.$
\end{lemma}

In the next Lemma we show identification of the states of $J$ and zeros of polynomial $F.$
\begin{lemma} \lb{L-ZerosF}
 The projection $\wt {}:\,\,\Lambda \mapsto \C$ of the set of states of  $J$ on $\Lambda$ coincides with the set of zeros
of $F$ on the complex plane $\C:$
$$\wt\s_{\rm st}\,(J)={\rm Zeros}\,(F).$$
Moreover, the multiplicities of  bound states and resonances are
equal to the multiplicities of zeros of $F.$ All bound states are
simple.
The virtual state is a simple zero of $F.$
\end{lemma}
{\bf Proof:}  First we observe that $f_0^+(\l)$ is analytic on
$\Lambda\setminus\s_{\rm st}\,(J^0)$.

 By Theorem\ref{th-charact-states} a
point $\l_0\in\g_k^+,$ $\l_0\not\in\sigma_{\rm st}(J^0),$ is a bound state
 iff $f_0^+(\l_0)=0.$
  Then $f_0^-(\l_0)\neq 0$ as the Wronskian
$\{f_0^+,f_0^-\}(\l_0)\neq 0.$
Moreover,  it follows that $\wt\l_0$ is zero
of $F(\l)$ with the same multiplicity (one).

A point $\l_0\in\Lambda_2,$ $\l_0\not\in\sigma_{\rm st}(J^0),$ $\Omega(\l_0)\neq 0,$
 is a resonance
 iff $f_0^+(\l_0)=0$ which  is equivalent to $f_0^-(\l_1)=0,$
 where $\l_1$ is the same number as $\l_0$ but on the physical sheet.
  Then it follows that $F(\wt\l_0)=0$ with the
same multiplicity.


 If $F(\l_0)=0$ for some
$\l_0\in\R,$ $\l_0\not\in\sigma_{\rm st}(J^0),$ $\Omega(\l_0)\neq 0,$ then
it is clear that there is either a bound state $\l_0^1\in\L_1$ with
$\wt\l_0^1=\l_0$ or an antibound  $\l_0^2\in\L_2$ state with $\wt\l_0^2=\l_0$
with the same multiplicity as $\l_0$.

If $F(\l_0)=0$ for some $\l_0\in\C\setminus\R,$  then necessarily
$f_0^+(\l_0^2)= 0$ at $\l_0^2\in\Lambda_2,$ with $\wt\l_0^2=\l_0,$
 and  $\l_0^2$ is the complex
resonance with the same multiplicity as $\l_0.$

Consider  now a point $\l_0\in\g_1^+$ or   $\l_0\in\g_1^-$ such that
 $\l_0\in\sigma_{\rm st}(J^0)$ and
  $\varphi^+_n(\wt\l_0)\neq 0$ for
some $n> 0.$  Then  $m_+$  has a pole at $\l_0,$ and
    $f_n^+(\l)$
 has a simple pole at $\l_0.$  Then $\l_0$ is a pole of
$${\mR}_n(\l)=\frac{f_n^+(\l)}{f_0^+(\l)}=
\frac{\vt^+_n+m_+\vp^+_n}{\vt^+_0+m_+\vp^+_0}$$
iff $\vp^+_0(\wt\l_0)=0,$ as by Lemma \ref{tph+} in this case  $\vt^+_0(\wt\l_0)\neq 0.$

Now using the identity $F_0=\vp_qf_0^+(\l)f_0^-(\l)=\vp_q(\vt^+_0)^2
+(\vp_{q+1}-\vt_q)
\vt^+_0\varphi^+_0-\vt_{q+1}(\varphi_0^+)^2$
 we get that if $\vp_q(\wt\l)=\vp^+_0(\wt\l)=0,$ then necessarily
 $ \wt\l$ is a simple zero of $F_0$ and $f_0^\pm(\l)\neq 0.$

The other statements of Lemma follows similarly as in the proof of Theorem \ref{th-charact-states} \hfill\BBox

Let $M_\pm\in\C$ denote (the projection of) the set of
poles  of $m_\pm.$ Let $M_{\rm e}$ denote the set of
 square root singularities of $m_\pm$ if $\mu_k=\l_k^+$ or $\mu_k=\l_k^-,$ $k=1,\ldots, q-1.$
 Note that $M_+\cap M_-=\emptyset.$
 We put
$$
{\displaystyle D^\pm=\prod_{\mu_k\in M_\pm} (\wt\l-\mu_k)},\qqq
{\displaystyle D^{\rm e}=\prod_{\mu_k\in M_{\rm e}}
\sqrt{\wt\l-\mu_k}},
$$ where $\,\wt {}\,:\,\,\Lambda\mapsto\C$ is the
natural projection introduced in (\ref{projection}). Let $\mu_\pm=\sharp\,( M_\pm),$ $\mu_{\rm
e}=\sharp\,( M_{\rm e}),$  be the number  of elements in the
respective sets.  If all gaps are open ($\l_n^-<\l_n^+,$
$n=1,\ldots,q$) then we have $\mu_++\mu_-+\mu_{\rm e}=q-1$ and
$\vp_q=a_0^0 (D^{\rm e})^2D^+D^-.$
We mark with $\hat{}$ the  modified (regularized) quantities:
$\hat{\psi}^\pm=D^{\rm e}D^\pm\psi^\pm,$ $\hat{f}^\pm = D^{\rm e} D^\pm f^\pm.$
Now $\hat{\psi}^\pm,$ $\hat{f}^\pm$
  are analytic in $\Lambda_1.$

In the next Lemma we prove the crucial property for the function
$F\equiv F_0=\vp_qf^+_0f^-_0=a_0^0\hat{f}^+_0\hat{f}^-_0.$  Recall that
$\{\phi_n,\psi_n\}={a}_n(\phi_n\psi_{n+1}-\phi_{n+1}\psi_n\}$
denotes the  Wronskian. Let as before $\dot y=\partial_\l y=\partial y/\partial \l$ and define the difference derivative
$$
\partial_n f(n)=f(n+1)-f(n).
$$

\begin{lemma}
\label{generalresults}
i) Any solution $y_n$ of (\ref{pert}) satisfies
\[\label{Fact1}
\pa_n\{\dot y,y\}_n=-( y_{n+1})^2,\,\,\forall n\ge 0.
\]
 ii) Suppose that $\l_1\in\g_k^+,$ for $k=0,1,\ldots,q$  and
$\hat{f}_0^+(\l_1)=0,$ i.e.  $\l_1$ is an eigenvalue of $J$ with the
eigenfunction  $y_n=\hat{f}_n^+(\l_1).$  Then
\begin{align}
&m_1:=\sum_{k=0}^\infty\left(\hat{f}_{k}^+(\lambda_1)\right)^2=
{a}_0\left(\frac{\partial}{\partial \l} \hat{f}_0^+\right)\hat{f}_1^+
>0\,\qq \mbox{at}\,\,\lambda=\lambda_1;
\label{Fact2}\\
&\{\hat{f}^+,\hat{f}^-\}_n=\vp_q(m_--m_+);
\label{Fact3}\\
&m_1=\frac{\dot{F}(\lambda_1)}{a_0^0(\hat{f}_0^-(\lambda_1))^2}
\cdot (-1)^{q-k+1}2\sinh
qh(\l_1)=\frac{(\partial_\l\hat{f}_0^+)(\l_1)}{\hat{f}_0^-(\lambda_1)}\cdot
(-1)^{q-k+1}2\sinh qh(\l_1)
>0\label{Fact6},
\end{align}
where $h(\l_1)=\Im \vk(\l_1)>0.$  Thus $(-1)^{q-k}\dot{F}(\l_1) <0$
and the function $F$ has simple zeros at all bound states of $J$ for
which $\varphi_q\neq 0.$ If $\l_0=\mu_k$ is an antibound state then
necessarily it is simple and $(-1)^{q-k}\dot{F}(\l_0)>0.$

\end{lemma}

{\bf Remark.} As by Lemma \ref{L-ZerosF} the zeros of $F$ coincide with the projections of states of $J$ to $\C$    then  Lemma \ref{generalresults} implies that between any
two  (projections of) eigenvalues $\l_{1}, \l_{3}\in\g_k$ (not separated by a band of the
absolute continuous spectrum) there is at least one (projection of) real resonance
(antibound state) $\l_2$ such that $(-1)^{q-k}\dot{F}(\l_2)>0.$

 {\bf Proof.}
 i)  Using
 $ y_{n+2}={1\/a_{n+1}}((\l-b_{n+1})y_{n+1}-a_{n}y_n),$ we get
  $$
  \pa_n\left[{a}_n(\dot
y_n) y_{n+1}-{a}_n(\dot y_{n+1})y_n\right]=-( y_{n+1})^2,
$$
 which yields (\ref{Fact1}).

ii) Note the following ``telescopic'' sum $\sum_{k=n}^m\partial
y_k=y_{m+1}-y_n.$ We put $n=0$ and get from (\ref{Fact1})
$$
\left\{\dot y,y\right\}_{m+1}-{a}_0\left[\left(\dot y_0\right)
y_1-\left(\dot y_1\right)y_0\right]=-\sum_{k=0}^m y_{k+1}^2.
$$
We put $\lambda=\lambda_1$ and $y=\hat{f}^+(\l_1).$
Then, using that the eigenfunction $ \hat{f}^+(\lambda_1)\in\ell^2(\N)$  and $\hat{f}_m^+\rightarrow 0$  as
$m\rightarrow\infty,$ we get that the first term in the left hand side goes to zero.
As  $\lambda_1$ is an eigenvalue, then  we have $\hat{f}_0^+(\lambda_1)=0$ and we
get
$$-{a}_0\left(\frac{\partial}{\partial
\l} \hat{f}_0^+\right)\hat{f}_1^+=-\sum_{k=0}^\infty
(\hat{f}_{k+1}^+)^2\,\,\mbox{at}\,\,\lambda=\lambda_1.$$  Finally
we get (\ref{Fact2})  using that $ \hat{f}^+(\lambda_1)\in\R.$

Next formula (\ref{Fact3})  follows from $\const=\{f^+_n,f^-_n\}=\{\psi^{+}_n,\psi^{-}_n\}=\{\psi^{+}_0,\psi^{-}_0\}={a}^0_0(m_--m_+).$

  Putting $n=0$  we get also
$\{f^+_n,f^-_n\}=-{a}_0f_1^+(\lambda_1)f_0^-(\lambda_1)$ using again
$f_0^+(\lambda_1)=0.$ Together with (\ref{Fact3}) and definitions of
$m_\pm$  it implies \begin{align}
&\hat{f}_1^+(\lambda_1)\hat{f}_0^-(\lambda_1)=
\frac{1}{a_0^0}\vp_qf_1^+(\lambda_1)f_0^-(\lambda_1)=\frac{\vp_q}{{a}_0}(m_+-m_-)=
\frac{i2\sin q\varkappa(\lambda_1)}{{a}_0}\nonumber\\
&\,\,\Rightarrow\,\,\hat{f}_1^+(\lambda_1)=\frac{i2\sin q\varkappa(\lambda_1)}{{a}_0\hat{f}_0^-(\lambda_1)}.\lb{Fact4}
\end{align}
Recall that $F(\lambda)=a_0^0\hat{f}_0^+\hat{f}_0^-.$ Taking the derivative of $F$  with respect to
$\l,$ we get $\dot{F}(\l_1)=a_0^0(\partial_\l\hat{f}_0^+)(\l_1)\hat{f}_0^-(\l_1),$ wherefrom it follows
\[\lb{Fact5}
(\partial_\l\hat{f}_0^+)(\lambda_1)=\frac{\dot{F}(\lambda_1)}{a_0^0\hat{f}_0^-(\lambda_1)}.
\]

Inserting (\ref{Fact4}) and (\ref{Fact5}) in (\ref{Fact2}):
 $m_1=\sum_{k=0}^\infty\left|\hat{f}_{k}^+(\lambda_1)\right|^2= {a}_0(\partial_\l
\hat{f}_0^+)(\l_1)\hat{f}_1^+(\l_1),$ we get
$$m_1=\dot{F}(\lambda_1)\cdot \frac{i2\sin q\varkappa(\lambda_1)  }{a_0^0(\hat{f}_0^-(\lambda_1))^2}
>0.$$

For $\lambda_1\in\g_k^+$ for $k=0,1\ldots,q,$
$\Im\varkappa(\lambda_1)= h(\l_1)>0.$  Then by (\ref{isin}) $i\sin
q\kappa(\l_1)=-(-1)^{q-k}\sinh qh(\l_1+i0),$ which implies
(\ref{Fact6}).

\hfill\BBox

\begin{lemma}\lb{l-general case}
i) The following identity holds true
\[
\lb{kvform} F=\vp_q\lt({\vt}_0^++ \frac{\f}{\vp_q}{\vp}_0^+\rt)^2
+{1-\D^2\/\vp_q}({\vp}_0^+)^2.
\]
Moreover, $F(\l)\neq 0,$ for any $\l\in (\l_{n-1}^+,\l_n^-),
n=1,\ldots,q,$ and $\sign F|_{(\l_{n-1}^+,\l_n^-)}=\sign\vp_q
|_{(\l_{n-1}^+,\l_n^-)}$.

ii) If $\l_0\in \{\l_{n-1}^+,\l_n^-\}$ is a virtual state, then $F$
has a simple zero at $\l_0$.

iii) There is always odd number $\ge 1$ of states (eigenvalues,
antibound or virtual state) in each finite open gap
$\g_n^c=\ol\g_n^-\cup \ol\g_n^+,$ $n=1,\ldots,q-1.$

\end{lemma}
{\bf Proof.} i) Using \er{TitchWeyl} and \er{mm} we obtain
$$
F=\vp_q\rt((\vt_0^+)^2+(m_++m_-)\vt_0^+\vp_0^++m_+m_-(\vp_0^+)^2\rt)
=\vp_q\rt((\vt_0^+)^2+{2\f\/\vp_q}\vt_0^+\vp_0^+-{\vt_{q+1}\/\vp_q}
(\vp_0^+)^2\rt)
$$
$$
=\vp_q\lt(\vt_0^++ \frac{\f}{\vp_q}\vp_0^+\rt)^2
+{\f^2-\vt_{q+1}\vp_q\/\vp_q}(\vp_0^+)^2 =\vp_q\lt(\vt_0^++
\frac{\f}{\vp_q}\vp_0^+\rt)^2 +{1-\D^2\/\vp_q}(\vp_0^+)^2.
$$
Now ii) and iii) follow directly from i).

 \hfill\BBox

 Now {\bf the proof of Theorem \ref{mainstatictheorem}} follows from the
properties of the function $F=\vp_qf^+f^-,$ stated in Lemmata \ref{defstates2}--\ref{l-general case}.

\hfill\BBox

In the next lemma we consider the zeros of the function $S(\l)-1$ which are the solutions of the equation $f_0^+(\l)=f_0^-(\l).$ Note that if $\l_1\in\L_1$ is a zero of $S-1$ then also $\l_2\in\L_2$ such that $\wt\l_2=\wt\l_1$ is a zero of $S-1.$

\begin{lemma}
\lb{S=1}Let $\l_0\in\L$ and $\wt\l_0\in\C$ denote the projection on $\L_1.$\\ i)  Suppose that $\vp_0^+(\wt\l_0)=0$  and one of the following conditions is satisfied:\\
$\phantom{==}$1) $\l_0\not\in\sigma_{\rm st}(J^0).$\\
$\phantom{==}$2) $\l_0\in\sigma_{\rm st}(J^0),$ $\Omega(\l_0)\neq 0$ and $\wt\l_0$ is zero of $\vp_0^+$ of multiplicity $\geq 2.$\\
$\phantom{==}$3) $\l_0\in\sigma_{\rm st}(J^0)$ and $\Omega(\l_0)=0.$

    Then
$S(\l_0)=1$.
\\ \\
ii) Suppose that $S(\l_0)=1$
and one of the following conditions is satisfied:\\
$\phantom{==}$1) $\l_0\not\in\sigma_{\rm st}(J^0)$ and $\Omega(\l_0)\neq 0.$   \\
$\phantom{==}$2)  $\l_0\in\sigma_{\rm st}(J^0)$ and $\Omega(\l_0)\neq 0.$ \\
$\phantom{==}$3) $\l_0\in\sigma_{\rm st}(J^0)$ and $\Omega(\l_0)= 0.$

 Then $\vp_0^+(\wt\l_0)=0.$\\ \\
 In the cases 1) and 3) the zeros (and their multiplicities) of $\vp_0^+$ and $1-S$ coincide.

\end{lemma}

{\bf Proof.} i)
Note the identities following from (\ref{P_2})
\begin{equation}\lb{smat}
1-S(\l_0)=\frac{f_0^+(\l_0)-f_0^-(\l_0)}{f_0^+(\l_0)}=\frac{2i\Omega(\l_0)}{\vp_q(\wt\l_0)}\frac{\vp_0^+(\wt\l_0)}{f_0^+(\l_0)}.
\end{equation}
Note that $\l_0\in\sigma_{\rm st}(J^0)$ iff $\vp_q(\wt\l_0)=0.$ Assume that $\vp_q(\wt\l_0)\ne 0.$  Then $f_0^\pm$ are
analytic at $\l_0$  and due to Lemma \ref{tph+} we obtain
$f_0^\pm(\l_0)=\vt_0^+(\wt\l_0)\neq 0$. Using this
 we
get $S(\l_0)={{f_0^-(\l_0)}\/f_0^+(\l_0)}=1$. This is also true for  $\Omega(\l_0)= 0.$

Assume now that $\l_0\in\sigma_{\rm st}(J^0).$   We distinguish between two cases.

Firstly, let $\wt\l_0\in \L_1$ be a zero of $\vp_0^+$ with multiplicity
$\ge 2$. Then $f_0^\pm(\l_0)=\vt_0^+(\wt\l_0)\ne 0$, since $\wt\l_0$ is a
simple zero of $\vp_q$. Thus
$S(\l_0)={{f_0^-(\l_0)}\/f_0^+(\l_0)}=1$.

Secondly, let $\wt\l_0\in \L_1$ be a simple zero of $\vp_0^+$. Suppose $\Omega(\l_0)\neq 0.$ As $\l_0\in\sigma_{\rm st}(J^0),$ then the point
$\l_0\in\L$ is a pole of $m_+$.
Then $m_-$ is analytic at $\l_0$ and using \er{TitchWeyl} we have
\begin{equation}\lb{dotphi}
f_0^+(\l_0)=\vt_0^+(\wt\l_0)+{2\f(\wt\l_0)\/\dot\vp_q(\wt\l_0)}\dot\vp_0^+(\wt\l_0),
\qqq f_0^-(\l_0)=\vt_0^+(\wt\l_0).
\end{equation}
This yields $f_0^+(\l_0)\ne f_0^-(\l_0)$, since
${2\f(\wt\l_0)\/\dot\vp_q(\wt\l_0)}\dot\vp_0^+(\wt\l_0)\ne 0$.   Note that
$\vt_0^+(\wt\l_0)\ne 0$. Then $S(\l_0)\ne 1$.

Suppose now that $\Omega(\l_0)= 0.$ Then $$m^\pm(\l)=\frac{c}{\sqrt{\epsilon}}+{\mathcal O}(1),\qqq \l=\l_0+\epsilon,\qqq\epsilon\rightarrow 0+,\qqq c\neq 0,$$
and $f^\pm_0(\l_0)=\vt_0^+(\wt\l_0)\neq 0$ which implies that $S(\l_0)=1.$

ii) Let $S(\l_0)= 1$. We use (\ref{P_2})
$$
\vp_0^+=\frac{\vp_q}{2i\O}\left(f_0^+-f_0^-\right)
=\frac{\vp_q}{2i\O}f_0^+\left(1-S\right).
$$
If $\Omega(\l_0)\neq 0$ and  $\vp_q(\wt\l_0)\neq 0,$ then $f^\pm_0$ are bounded near $\l_0$  and we have $\vp_0^+(\wt\l_0)=0.$

If $\Omega(\l_0)\neq 0$ and $\l_0\in\sigma_{\rm st}(J^0),$  then $\wt\l_0$ is the zero of $\vp_0^+,$ and from (\ref{dotphi}) it follows that the multiplicity  of $\wt\l_0$ is $\geq 2.$

 If $\Omega(\l_0)=0$ and $\vp_q(\wt\l_0)\neq 0,$ then $\vp_0^+(\wt\l_0)\neq 0.$

 If $\Omega(\l_0)=0$ and $\l_0\in\sigma_{\rm st}(J^0),$  then  we get $f_0^+(\l_0)=\vt_0^+(\wt\l_0)\neq 0$ and $\vp_0^+(\wt\l_0)= 0$ as $\Omega(\l)=c\sqrt{\epsilon}+{\mathcal O}(\epsilon)$ as $\l-\l_0=\epsilon\rightarrow 0+.$
\hfill\BBox

\section{Inverse problem}\label{s-inverse}
\setcounter{equation}{0}

\subsection{Preliminaries}\lb{ss-inverse-prel}
In this section we collect some properties of the Jost solutions needed for the proof of the inverse results.
The first lemma states that that the Jost solutions $f^\pm$  inherit the properties of $\p^\pm.$
\begin{lemma}
\lb{l-jost} 1) Each $f_n^\pm, n\ge 0$, is analytic in $\cZ\sm \{0\}$
and continuous up to $\partial
\cZ\setminus\{z(\mu_j)\}_{j=1}^{q-1}.$ Moreover, the following
identities hold true:
\[
f^\sigma=\vt^\sigma +m_\sigma \vp^\sigma, \qqq \sigma=\pm.
\]
\[
 f_n^\pm(\overline{z})=f_n^\pm(z^{-1})=f_n^\mp(z)=\overline{f_n^{\, \pm}(z)}\qq \mbox{for}\qqq |z| =1.
\]
2) $f_n^\pm(z)$ does not have a singularity at $z (\mu_j)$  if $\mu_j$ is not a singularity
(square root singularity if $\mu_j$ coincides with the band edge) of  $m_{\pm},$ otherwise, $f_n^\pm(z)$  can have either a simple pole at $z (\mu_j)$ if $\mu_j$ is a pole of $m_\pm,$ or a square root singularity,
\[
\lb{B2} f_n^\pm(\l)=\pm \sigma(-1)^{q-j}\frac{iC(n)}{\sqrt{\l -\l_j^\sigma}} +\cO(1),\qq \l\in [\l_{j-1}^+,\l_j^-],
\]
if $\mu_j$  coincides with the band edge:  $\mu_j=\l_j^\sigma,$ $\sigma=+$ or $\sigma=-$, $j=1,\ldots, q-1.$ Here
  $C(n)$ is bounded and real, the factor $ \sigma(-1)^{q-j}$ comes from the analytic continuation of the square root $\Omega(\l)$ using Definition
(\ref{branch}).
\end{lemma}

The asymptotics of the function $f^+(z)$ are given in
(\ref{asjo1}), (\ref{asjo2}).

The next lemma follows from the straightforward reformulation of the results obtained in Section \ref{s-genJost} in the form stated in the definition of  $\gJ_\nu.$
\begin{lemma} If $(u,v)\in\gX_\nu,$  where $\nu=2p$ or $\nu=2p-1,$
 then the Jost functions $f_0^\pm\in \gJ_\nu$ (see Definition \ref{classJ}).
\end{lemma}

\subsection{ Inverse scattering problem.}\lb{s-scatt}  In this section we recall
some relevant for us results from \cite{Kh2} and \cite{EMT}. Let
$\hat{S}=\frac{\hat{f}^-(\l)}{\hat{f}^+(\l)}.$ Then the scattering
matrix is $S=\frac{D^+}{D^-}\hat{S}.$ For each eigenvalue $\gr_n$ we
define the norming constant $m_n$ by
\[
m_n=\sum_{j=0}^\iy\left(\hat{f}_j^+(\gr_n)\right)^2, \qqq n=1,\ldots,N.
\]

Introduce the scattering data for the pair of operators $J, J^0$ by
$$
\cS(J)=\rt\{\hat{S}(\l),\,\,\mbox{for}\,\,\l\in\s_{\rm
ac}(J^0),\,\,\gr_k,\,\,m_k,\,\,k=1,2,\ldots,N\rt\}.
$$
 By the inverse scattering theory for this pair, we understand  the problem of
reconstructing the perturbed operator $J$ from the scattering data
and the unperturbed operator $J^0.$

Everywhere in this section we assumes  that $(u,v)\in\gX_\n.$ We
introduce the Gel'fand-Levitan-Marchenko  equation for a matrix $K(n,m)$  by
\[
\lb{GLM2} K(n,m)+\sum_{l=n}^{+\infty} K(n,l)
\gF_{l,m}=\frac{\d_{nm}}{K(n,n)},\qq  m\geq n.
\]
Here the sum in (\ref{GLM2}) is finite, since $(u,v)\in\gX_\n$. The
matrix $\gF_{l,m}$ is constructed from the  scattering data
$\cS(J)$ by
\[
\lb{7.7}
\gF_{l,m}=\gF_{l,m}^0+\sum_{j=1}^N{\hat\p^+_l(\gr_j)\hat\p_m^+(\gr_j)\/m_j}
,
\]
where
$$
\gF_{l,m}^0=-\frac{1}{2\pi i}\int_{|z|=1} S(z) \psi^+_l(z)
\psi_m^+(z)d\omega(z)
$$
and
\[\lb{domega}
d\omega(z)=\prod_{j=1}^{q-1}\frac{\l(z)-\mu_j}{\l(z)-\alpha_j}\frac{dz}{z}.
\]
Here $\alpha_j\in\g_j^+$ is the zero of $\Delta'(\l)$
 (see Section \ref{ss-cut} and  (3.22) in \cite{EMT} ).
Note that $\gF_{l,m}^0=\gF_{m,l}^0$ and $\gF_{l,m}^0$ is real. We
will determine the  matrix $K(n,m)$  from the Gel'fand-Levitan-Marchenko equation
\er{GLM2} and reconstruct (see (5.27) in \cite{EMT})  $a_n, b_n$ by
\[
\label{5.27}
\frac{{a}_n}{a^0_n}=\frac{K(n+1,n+1)}{K(n,n)},\qq b_n= b_n^0+
a^0_n\frac{K(n,n+1)}{K(n,n)}-a^0_{n-1}\frac{K(n-1,n)}{K(n-1,n-1)}.
\]

Now we  consider the Gel'fand-Levitan-Marchenko equation. From \cite{Kh1} or
\cite{EMT}, Lemma 5.1, it is known that the Jost solution $f^+_n$
can be represented as
$$f_n^+(z)=\sum_{m=n}^\infty K(n,m)\psi_m^+(z),\,\,|z|=1,$$
where for $(u,v)\in\gX_\n$ the kernel $K(n,m)$ has finite rank and
satisfies
$$
K(n,m)=0,\qqq {\rm for} \qqq m<n,
$$
\[\label{5.5}
|K(n,m)|\leq C\sum_{j=[\frac{m+n}{2}]+1}^p\left(|{u} _j|+|{v} _j|\right),\,\,m>n.
\]
Here the constant $C\equiv C(J^0)$ depends on the unperturbed operator $J^0.$

We recall the properties of the scattering data $\cS(J)$ from
\cite{Kh2}.

\begin{enumerate}
\item[(I)]
{\em Function $S(\l)$ is continuous for $\l\in{\rm int}\, \partial\G,$ where $\G$ is the cut plane $\C\setminus\sigma_{\rm ac}(J^0),$
 $$
 \overline{S(\l)}=S^{-1}(\l),\,\,\l\in {\rm int}\,\partial\G,\,\,
 \mbox{and}\,\,S(\l-i0)=\overline{S(\l+i0)},\,\,\l\in{\rm int}\,
 \s_{\rm ac}(J^0),
 $$
 where
 ${\rm int}$ stands for interior.}
\item[(II)] {\em The function
$$\gF_{l,m}^0=-{1\/2\pi i}\int_{|z|=1}S(z) \p^+_l(z) \p_m^+(z)d\o(z)
$$
satisfies
\[\lb{sum}
\sum_{l=0}^\iy\sup_{m\geq 0}|\gF_{l,m}^0|<\iy.
\] In \cite{Kh2} this function was denoted $S(n,m).$}
\item[(III)] {\em  Equation
\[
\lb{han15}
h_m+\sum_{k=1}^\infty S_{m,k}h_k=0,\,\,m=1,2,\ldots,
 \]
 has precisely $N$ linearly independent solutions in $\ell^2(1,\infty).$}
\item[(IV)] {\em The equation $\sum_{m=-\infty}^0S_{l,m}g_m=g_n$ has only
the zero solution in $\ell^2(-\infty,0).$}
\item[(V)]
{\em The quantities $a_n$ and $b_n$ defined in (\ref{5.27}),
where $K(n,m)$ is solution to (\ref{GLM2}), satisfy the inequality
$$
\sum_{n=1}^\infty n\left( \left|\left(\frac{{a}_n}{a_n^0}\right)^2-1\right|+|b_n-b_n^0|\right) <\infty.
$$}
\end{enumerate}

\begin{theorem}[Khanmamedov]\lb{ThHan}
If conditions (I)--(III) hold, then for every $n\in\N,$
the Gel'fand-Levitan-Marchenko equation (\ref{GLM2}) has unique solution in
$\ell^2(n+1,\infty).$

The set $\cS(J)$  uniquely determines $J$ iff conditions (I)--(V) hold.
\end{theorem}

From the proof of Khanmamedov it follows that:\\ if $({u} ,{v}
)\in\gX_\n,$ the bound states $\gr_j\in\g_k,$ $k=0,\ldots,q,$  the
norming constants $m_k$ are given by
$m_j=\sum_{n=0}^\infty\left(\hat{f}_{n}^+(\gr_j)\right)^2$ and
$S$--matrix is  given by $S=\frac{f^-_0(\l)}{f^+_0(\l)},$\\ then
conditions (I)--(V) are satisfied.

Recall that from Lemma \ref{generalresults}, property (\ref{Fact6}), it follows that for
$\gr_j\in\s_{\rm bc}\cap\g_k^+$ we have
\[\lb{norming}
m_j=\frac{\dot{F}(\gr_j)}{a_0^0(\hat{f}_0^-(\gr_j))^2}\cdot
(-1)^{q-k+1}2\sinh
qh(\gr_j)=\frac{(\partial_\l\hat{f}_0^+)(\gr_j)}{\hat{f}_0^-(\gr_j)}
(-1)^{q-k+1}2\sinh qh(\gr_j)
>0,
\]
 where $h(\gr_j)=\Im\vk(\gr_j) >0$ (see (\ref{isin})),
 as $\dot{F}(\gr_j)=a_0^0(\partial_\l\hat{f}_0^+)(\gr_j)\hat{f}_0^-(\gr_j),$ $(-1)^{q-k}\dot{F}(\gr_j)=a_0^0(-1)^{q-k}(\partial_\l\hat{f}_0^+)(\gr_j)\hat{f}_0^-(\gr_j) <0.$

Now we show that the scattering data $\cS$ can be uniquely reconstructed from any function  $f\in \gJ_\nu$
as in  Definition \ref{classJ} and the
conditions
(I)--(V) are satisfied.

 The $S-$matrix and the norming constants $m_j,$ $1\leq j\leq N,$ are then expressed in terms of the function  $f\in \gJ_\nu$ only. By abuse of notation we will keep the same letters $S$ and $m_j$ for the functions  expressed in $f.$

Using Theorem \ref{ThHan} this  will  imply that the function $f\in \gJ_\nu$ uniquely determines $J.$

\begin{lemma}\lb{LthHan}  Let  $f=P_1+m_+P_2\in\gJ_\nu,$ $f_-=P_1+m_-P_2,$
$P(\l)=\vp_q ff_-$ and  $\s_{\rm bs}(f)=\{\gr_j\}_{j=1}^N\in\L_1$ be
as in Definition \ref{classJ}. We define $m_j,$ $j=1,\ldots,N,$ by
\begin{equation}
\lb{norming-f} m_j=\frac{\dot{P}
(\gr_j)}{a_0^0(\hat{f}_-(\gr_j))^2}\cdot
(-1)^{q-k+1}2\sinh
qh(\gr_j),\qq\qq\gr_j\in\gamma_k^+,\end{equation}
where $\hat{f}_-=D^{\rm e}D^-f_-,$
 and
$S(\l):=\frac{f_-(\l)}{f(\l)}.$
    Then conditions (I)-(V) are satisfied.
\end{lemma}
{\bf Proof.} (I) Recall that by  \er{smatrix}
$S(\l)={\ol{f_0^+(\l)}\/f_0^+(\l)}={{f_0^-(\l)}\/f_0^+(\l)}$, and
then it follows
$$\overline{S(\l)}=S^{-1}(\l),\,\,\l\in {\rm
int}\,\partial\G,\,\,\mbox{and}\,\,S(\l-i0)=\overline{S(\l+i0)},\,\,\l\in{\rm
int}\,\s_{\rm ac}(J^0),
$$

(II)  In the next section we prove  that if
$\{\l_j\}_{j=1}^\k\in\s_{\rm st}(f)$ then the sum (\ref{sum}) is finite
and the condition is trivially satisfied.

(III)  Khanmamedov \cite{Kh2} showed that the number of linearly
independent solutions in $\ell^2(1,\infty)$ of (\ref{han15})
coincides with that of linearly independent functions of the form
$\frac{C_k\hat{f}_0^+(\l)}{\partial_\l\hat{f}_0^+(\gr_j)(\l-\gr_j)}.$
For $\{\l_j\}_{j=1}^\k\in\s_{\rm st}(f)$ as in Introduction it follows
that  the values $\gr_j\in\R\setminus\s_{\rm ac}(J^0),$ $1\leq j\leq
N,$ are distinct and the norming constants $m_j,$ $1\leq j\leq N,$
are positive, which implies that the number of linearly independent
functions is precisely $N.$

(IV) The condition is proved similarly to (III).

(V) For $({u} ,{v} )\in\gX_\nu$ and $a_n,$ $b_n$ defined in (\ref{ab})
or
for  $\{\l_j\}_{j=1}^\k\in\s_{\rm st}(f)$ for $f\in\gJ_\nu,$ as in Definition \ref{classJ}, this sum
is finite as shown in the next section. \hfill\BBox

\

\subsection{Inverse resonance problem. }\lb{s-invres} We prove here
the Theorems \ref{th-inverse}-\ref{th_getJost_2}.

{\bf Proof of Theorem \ref{th-inverse}.}

We will prove the following: \no {\it  The mapping
$\mF  : \gX_\nu \to \gJ_\nu$ given by
$$
 ({u} ,{v} )\to f_0^+({u} ,{v} )\in \gJ_\nu,
 $$
 is
one-to-one and onto. Recall that $\nu\in\{2p-1, 2p\}$.  In
particular,  a pair of coefficients in $\gX_\nu$ is uniquely
determined by its bound states and resonances.}

{\bf Uniqueness.} In the first part of this paper we proved  that to
any $({u} ,{v} )\in\gX_\n$ we can associate the Jost function
$f\in\gJ_\nu.$ Let $\s_{\rm st}(f)$ be the class of points on
$\L$ specified in Definition \ref{classJ},   $f_-=P_1+m_-P_2,$
 the bound states $\gr_j\in\s_{\rm
bs}(f)\subset\Lambda_1,$  the norming constants $m_j$ by
(\ref{norming-f}),  $j=1,\ldots,N,$ and the scattering matrix
$S=\frac{f_-}{f}.$ Then conditions (I)--(V) of Theorem
\ref{ThHan}  are satisfied and these data determine $({u} ,{v}
)\in\gX_\n$ uniquely. Then we have that the mapping  $({u} ,{v} )\to
f_0^+({u} ,{v} )\in \gJ_\nu$ is an injection.

{\bf Surjection.} We will show  that the mapping $ ({u} ,{v} )\to f_0^+({u} ,{v} )\in \gJ_\nu$ is surjective. Let
$f\in\gJ_\nu$ as in Definition \ref{classJ}.

Then we define $m_j,$ $j=1,\ldots,N,$ by (\ref{norming-f}) and
$\hat{S}=\frac{\hat{f}_-}{\hat{f}},$ where $\hat{f}=D^+D^{\rm e}f,$ $\hat{f}_-=D^-D^{\rm e}f_-.$  Lemma \ref{LthHan} shows
that  the set of quantities
$\cS=\{\hat{S}(\l),\,\,\mbox{for}\,\,\l\in\s_{\rm ac}(f),\,\,z_k,\,\,
m_k,\,\,k=1,2,\ldots,N\}$ is unique scattering data verifying
conditions (I)--(V). Then by solving the Gel'fand-Levitan-Marchenko equation and
applying  Theorem \ref{ThHan} we get the unique coefficients $({u} ,{v}
).$    We need to show that $({u} ,{v} )\in \gX_\nu.$

We have
 \begin{align*}
 \gF^0_{l,m}&=-\frac{1}{2\pi i}\int_{|z|=1}
S(z) \psi^+_l(z) \psi_m^+(z)d\omega(z),\\
&=-\frac{1}{2\pi i}\int_{\partial\G} \hat{S}(\l) \hat{\psi}^+_l(\l)
\hat{\psi}_m^+(\l) \frac{d\l}{2(\D^2(\l)-1)^{1/2}},
\end{align*}

Observe that $d\omega$ is meromorphic on $\cZ_1$ with simple pole at $z=0.$
 In particular, there are no poles at $z(\alpha_j).$
To evaluate the integral we use the residue theorem. Take a closed
contour in $\cZ_1$ and let this contour approach $\partial \cZ_1.$
The function $S(z) \psi^\pm_l(z) \psi_m^\pm(z)$ is continuous on
$\{|z|=1\}\setminus\{z(E_j)\}$ and meromorphic on $\cZ_1$ with
simple poles at $z(\gr_j)$ and eventually a pole at $z=0.$

We have
$$ S(z)=z^{-\nu}\left(1+{\mathcal O}(z)\right),\qqq\psi_l^+\psi_m^+=z^{l+m}\left(1+{\mathcal O}(z)\right),\qqq\mbox{as}\qq z\rightarrow 0.$$ Suppose $l+m\geq \nu+1$  ($+1$ is due to singularity of $z^{-1}$ in $d\omega$). Then the integrand is bounded near $z=0$ and we
apply the residue theorem to the only poles at the eigenvalues.

We have (\cite{EMT}, (3.23))
$$\frac{dz}{d\l}=z\frac{\prod_{j=1}^{q-1}(\l-\alpha_j)}{2(\Delta^2(\l)-1)^{1/2}}$$
and, if $z_j=z(\gr_j),$ then ${\rm Res}_{z=z_j} F(z)=z'(\gr_j){\rm
Res}_{\l=\gr_j}F(z(\l)).$

We get
$$
\gF_{l,m}^0=-\sum_{j=1}^N{\rm Res}_{\gr_j}\left(\frac{\hat{S}(\l)
\hat{\psi}^+_l(\l)
\hat{\psi}_m^+(\l)}{2(\Delta^2(\l)-1)^{1/2}}\right),
$$
where
$(\Delta^2(\l)-1)^{1/2}=i\Omega(\l).$ Now
$$
\hat{S}(\l)=
\frac{\hat{f}_-(\gr_j)}{\partial_\l\hat{f}(\gr_j)(\l-\gr_j)}\left(1+{\mathcal O}(\l-\gr_j)\right)\qqq\mbox{as}\qq \l\rightarrow \gr_j.
$$
Then, using that $2(\Delta^2(\l)-1)^{1/2}= (-1)^{q-k+1}2\sinh
qh(\l)$ (see (\ref{isin})) and (\ref{norming-f}), we get
\begin{align*}
\gF_{l,m}^0= -\sum_{j=1}^N\frac{\hat{f}_-(\gr_j)}{\pa_\l\hat{f}
(\gr_j)2(\Delta^2(\gr_j)-1)^{1/2}} \hat{\psi}_l^+(\gr_j) \hat{\psi}_m^+(\gr_j)
=-\sum_{j=1}^Nm_j^{-1}\hat{\psi}_l^+(\gr_j) \hat{\psi}_m^+(\gr_j)
\end{align*}
Then equation (\ref{7.7}) implies
\begin{equation*}
\gF_{l,m}=\gF_{l,m}^0+\sum_{j=1}^Nm_j^{-1}\hat{\psi}^+_l(\gr_j)
\hat{\psi}_m^+(\gr_j)=0,\qq l+m\geq \nu+1,
\end{equation*}
and the
Gel'fand-Levitan-Marchenko equation
$$
K(n,m)+\sum_{l=n}^{+\infty} K (n,l)
\gF_{l,m}=\frac{\delta_{nm}}{K(n,n)},\qq  m\geq n,
$$
implies that the  kernel of the transformation operator
$K(n,m),$ satisfies
$$
K(n,m) =\frac{\delta_{nm}}{K(n,n)},\qq  m\geq n,\qq m+n\geq \nu+1.
$$
Thus we get
$$
\mbox{If}\,\,2n\geq \nu+1,\,\,\mbox{then}\,\,K(n,n)=\pm 1;\qq
\mbox{if}\,\,n+m\geq \nu+1,\,\,m\neq n,\,\,\mbox{then}\,\,K(n,m)=0.
$$

We recall (\ref{5.27})
$$\frac{{a}_n}{a^0_n}=\frac{K(n+1,n+1)}{K(n,n)},\qq{v}_n=
a^0_n\frac{K(n,n+1)}{K(n,n)}-a^0_{n-1}\frac{K(n-1,n)}{K(n-1,n-1)}.
$$
Then, as ${a}_n > 0,$ ${a}^0_n>0,$  we get ${a}_n={a}^0_n$  for $n\geq p+1,$ if $\nu=2p$ (or for $n\geq p$ if $\nu=2p-1$). Moreover, we get  ${v}_n=0$ for $2n-1\geq 2p+1$ (or $2n-1\geq 2p$) which both implies $n\geq p+1$ and ${v} _p\neq 0,$ if $\nu=2p-1.$  This yields surjection.

From (\ref{5.5}) we get also
 that if $({u} ,{v} )\in \gX_\nu$ then $K(n,m)=0$ for $n+m\geq 2p.$
\hfill\BBox

{\bf Proof of Theorems \ref{th_getJost_1} and \ref{th_getJost_2}.} Recall that for any $\l\in\L$ the map $\l\mapsto\wt\l\in\L_1$ denotes the projection to the first sheet and $\L_1$ is identified with $\Gamma=\C\setminus\sigma_{\rm ac}(J^0).$
Note that from Lemma \ref{S=1} it follows that, due to  the assumption (\ref{hypS}), the (projection of) set of solutions of the equation $S=1$ is the set of all zeros of polynomial $\vp_0^+.$
 Recall that the polynomials $F,$ $\vp_0^+$ have orders $\k=\nu+q-1$ and $\nu-1,$ respectively.  We denote their sets of zeros by ${\rm Zeros}\,(F)=\{\l_j\}$ and ${\rm Zeros}\,(\vp_0^+)=\{\omega_j\}$ respectively.
 Now for given ${\rm Zeros}\,(F),$  ${\rm Zeros}\,(\vp_0^+)$ and the constants $c_1,$ $c_2,$ we can reconstruct the unique polynomials $F(\l)=C_1\prod_{j=1}^{\nu+q-1}(\l-\l_j),$ $\vp_0^+(\l)=C_2\prod_{j=1}^{\nu-1}(\l-\omega_j).$ We need to distinguish between projections to the complex plane  of the bound states and the resonances.

 Let $\sigma_1=\{\l_j\}_{j=1}^{N_1},$ $N_1\leq N,$ be the set of zeros of $F$   such that:\\
 1)  $\s_1\cap\wt\s_{\rm bs}(J_0)=\emptyset;$\\
 2) $\sigma_1\in\bigcup_0^q\gamma_j $  and if
  $\l_0\in\wt\sigma_1\bigcap\gamma_n$ for some $n=0,\ldots q,$ then  $(-1)^{q-n}\dot{F}(\l_0)<0.$

  Let $\sigma_2=\{\l_j\}_{j=N_1+1}^{\k_1},$ $\k_1\leq\k,$ be the set of zeros of $F$ such that:\\
   1) $\s_2\cap(\wt\s_{\rm r}(J_0)\cup \wt\s_{\rm vs}(J_0))=\emptyset;$\\
    2)  if  $\l_j\in\sigma_2$ is real, then  $\l_j\in\gamma_n,$ for some $n=0,\ldots q,$ and $(-1)^{q-n}\dot{F}(\l_j)\geq 0.$

  We consider the following polynomial interpolation problem:
    \begin{align}
     &\vt_0^+(\l_j)=-m_+(\l_j)\vp_0^+(\l_j)\qq \mbox{for}\qq \l_j\in\s_1,\,\, j=1,\ldots,N_1,\nonumber\\
    &\vt_0^+(\l_j)=-m_-(\l_j)\vp_0^+(\l_j)\qq \mbox{for}\qq \l_j\in\s_2,\,\, j=N_1+1,\ldots,\k_1.\lb{polint}
    \end{align}
    Suppose that each zero $\l_j\in\sigma_1\cup\sigma_2\subset{\rm Zeros}\,(F),$ $j=1,\ldots,\k_1,$ is simple. Then we have $\nu\leq \k_1\leq \nu +q-1$ and it is well known (see for example the book of Kendell A. Atkinson \cite{A}) that the polynomial interpolations problem (\ref{polint}) defines unique polynomial $\vt_0^+$ of order $\nu-2,$ and therefore the unique Jost function $f_0^+=\vt_0^++m_+\vp_0^+.$
    \hfill\BBox


\section{Asymptotics of the Jost function on the unphysical sheet.}\label{s-asym}

In this section  we obtain asymptotics of the Jost solutions $f^\pm$ and prove  Lemma \ref{Tm}. The asymptotics of $f^+(\l)$ as $\l\in\Lambda_1$ and $\l\rightarrow\infty$ are well known (see for example \cite{T}). We obtain the asymptotics of $f_{p-n}^+(\l)$ as $\l\in\L_2$ and
$\l\rightarrow\infty,$ which is equivalent to the asymptotics of
$f_{p-n}^-$ for $\l\in\Lambda_1.$ In this section we will not assume
$A=1.$  We will omit the upper indexes $\phantom{=}^\pm$ as much as
possible.  We make use of (\ref{defJost}):
 $$
 f_{p+1}=\psi_{p+1},\qq
f_p=\frac{a^0_p}{{a}_p}\psi_p.
$$
Put
$\Phi(j)=\frac{\psi_{j+1}}{\psi_j}.$ Now (see \cite{T}) we have
$$
\psi_p={\prod_{j=0}^{p-1}}{}^*\Phi(j)=
\left\{
                                          \begin{array}{ll}
                                            {\prod_{j=0}^{p-1}}{}\Phi(j) & \mbox{for}\,\,p>0 \\
                                            1 & \mbox{for}\,\, p=0 \\
                                            {\prod_{j=0}^{p-1}}{}(\Phi(j))^{-1}. &
\mbox{for}\,\, p<0,
                                          \end{array}
                                        \right.
$$
If $\psi=\psi^\pm$  then $\Phi(0)=\Phi^\pm(0)=m_\pm$ and we have (see \cite{T})
$$
\Phi^\pm(\l,n)=\left(\frac{a^0(n)}{\l}\right)^{\pm
1}\left(1\pm\frac{b^0(n+\!\!{\scriptsize
                                               \begin{array}{c}
                                                 1 \\
                                                 0 \\
                                               \end{array}}\!\!)}{\l}
+\cO\left(\frac{1}{\l^2}\right)
 \right),\qqq\l\rightarrow\infty,
                                             $$
where $a_n^0\equiv a^0(n),$ $b_n^0\equiv b^0(n).$
Put $\Psi(n)=\Phi^{-1}(n),$ then
$$
\Psi^\pm(\l,n)=\left(\frac{a^0(n)}{\l}\right)^{\mp
1}\left(1\mp\frac{b^0(n+\!\!{\scriptsize
                                               \begin{array}{c}
                                                 1 \\
                                                 0 \\
                                               \end{array}}\!\!)}{\l}
+\cO\left(\frac{1}{\l^2}\right)
 \right),\qqq\l\rightarrow\infty.
                                             $$
By iterating the Jacobi equation (\ref{eq-pert}) we get
\begin{align*}
f_{p-1}&=\frac{(\l-{b}_p)a^0_p\psi_p-{a}_p^2\psi_{p+1}}{{a}_p{a}_{p-1}}=
\frac{\psi_{p+1}}{{a}_p{a}_{p-1}}\left(
(\l-{b}_p)a^0_p\Psi(p)-{a}_p^2\right);\\
f_{p-2}&=\frac{(\l-{b}_{p-1}){a}_{p-1}f_{p-1}-
{a}_{p-1}^2\frac{a^0_p}{{a}_p}\psi_{p}}{{a}_{p-1}{a}_{p-2}}=\\
&=\frac{\psi_{p+1}}{{a}_p{a}_{p-1}{a}_{p-2}}\left((\l-{b}_{p-1})\left[(\l-{b}_p)a^0_p\Psi(p)-{a}_p^2\right]
-{a}^2_{p-1}{a}^0_p\Psi(p)\right);\\
f_{p-3}&=\frac{(\l-{b}_{p-2}){a}_{p-2} f_{p-2}-
{a}_{p-2}^2\frac{\psi_{p+1}}{{a}_p{a}_{p-1}}\left(
(\l-{b}_p)a^0_p\Psi(p)-{a}_p^2\right)}{{a}_{p-2}{a}_{p-3}}=\\
&=\frac{\psi_{p+1}}{{a}_p\ldots{a}_{p-3}}\left(
(\l-{b}_{p-2})\left[(\l-{b}_{p-1})\left[(\l-{b}_p)a^0_p\Psi(p)-{a}_p^2\right]
-{a}^2_{p-1}{a}^0_p\Psi(p)\right]-\right.\\
&-\left.{a}_{p-2}^2\left(
(\l-{b}_p)a^0_p\Psi(p)-{a}_p^2\right)\right).
\end{align*}
Now we use that $\displaystyle
\Psi(p)\equiv\Psi^-(\l,p)=\frac{a^0_p}{\l}\left(1+\frac{b^0_p}{\l}
+\cO\left(\frac{1}{\l^2}\right)
 \right)$ as $\l\rightarrow\infty.$
 Then we get
$$
\psi_{p+1}\equiv\psi_{p+1}^-(\l)=\frac{\l^{p+1}}{A_p}
\left(1-\frac{1}{\l}\sum_{j=0}^{p}b^0_j
+\cO\left(\frac{1}{\l^2}\right)
 \right),\qqq\l\rightarrow\infty,
$$
where $A_p=\prod_{j=0}^{p}a^0_j$. \ We have
$$
(\l-{b}_p)a^0_p\Psi(p)-{a}_p^2=((a^0_p)^2-{a}_p^2)+
\frac{(a^0_p)^2}{\l}(b^0_p-{b}_p)+\cO\left(\frac{1}{\l^2}\right)
$$ and get
\begin{align}
f_{p-n} &=\frac{\l^{p+n}}{ A_p\prod_{j=p-n}^p{a}_j}\nonumber\\
&\cdot\left(((a^0_p)^2-{a}_p^2)+{1\/\l}\left[
-((a^0_p)^2-{a}_p^2)(\sum_{j=0}^{p}b^0_j+\sum_{j=p-n+1}^{p-1}{b}_j)-
(a^0_p)^2v_p\right]+{\cO(1)\/\l^2}\right),
\nonumber\\
f_0(\l)&=\frac{c_1\l^{2p}}{A_p}\label{asJost}\\
&\cdot\left( ((a_p^0)^2-{a}_p^2)+\l^{-1}\left[
-((a^0_p)^2-{a}_p^2)(\sum_{j=0}^{p}b^0_j+\sum_{j=1}^{p-1}{b}_j)-
(a^0_p)^2v_p\right]+{\cO(1)\/\l^{2}}\right).\nonumber
\end{align}
If ${a}_p={a}^0_p,$ then
$$
f_0(\l)=-\frac{c_1(a_0^p)^2v_p}{{A_p}}\l^{2p-1}+{\mathcal
O}\left(\l^{2p-2}\right).
$$

Multiplying
\begin{align*}
\vp_q&=\frac{\l^{q-1}}{A_{q-1}}+\cO(\l^{q-2}),\\
f_n^+&=\alpha_n^+\frac{\prod_{j=0}^{n-1}{}^*a_j}{\l^n}
\left[1+\frac{1}{\l}\left(-\sum_{j=1}^p{v}_j+
\sum_{j=1}^{n}{}^*b_j\right)
+{\cO(1)\/\l^{2}}\right],\\
(f_n^+)^*&=\frac{\l^{2p-n}}{\prod_{j=n}^p{a}_jA_p}\\
&\cdot\left(  ((a^0_p)^2-{a}_p^2)+\l^{-1}\left[
({a}_p^2-(a^0_p)^2)(\sum_{j=0}^{p}b^0_j+\sum_{j=n+1}^{p-1}{b}_j)-
(a^0_p)^2v_p\right]+{\cO(1)\/\l^{2}}\right),
 \end{align*}
and using $\displaystyle
\alpha_n^+=\prod_{j=n}^p\frac{a^0_j}{{a}_j},$ we get
\begin{align}
F_n(\l)&=\vp_qf_n^+(f_n^+)^*= \frac{c_1^2\l^{2(p-n)+q-1}}{A_{q-1}
}\left(((a^0_p)^2-{a}_p^2)+{\mathcal O}(\l^{-1})\right), \qq \mbox{if} \qq
u_p\ne 0, \lb{F_n_1}
\end{align}
\begin{align}
F_n(\l)&=\vp_qf_n^+(f_n^+)^*= \frac{c_1^2\l^{2(p-n)+q-2}}{A_{q-1}
}\left(-(a^0_p)^2{v}_p+{\mathcal O}(\l^{-1})\right), \qq \mbox{if} \qq
u_p=0, \qq v_p\ne 0,\lb{F_n_2}
\end{align}
where $c_1(n)=(\prod_{j=n}^p{a}_j)^{-1}.$
On the Riemann surface $\cZ$ as in Section \ref{ss-inverse-prel} we get
\begin{align}
f_0^+= &\alpha_0^++\cO(z),\,\,\mbox{as}\,\, z\rightarrow 0,\label{asjo1}\\
f^+_0=&\frac{c_1A^{\frac{2p}{q}}z^{2p}}{A_p}\label{asjo2}\\
&\cdot\left(  ((a^0_p)^2-{a}_p^2)+{A^{-{1\/q}}\/z}\left[
-((a^0_p)^2-{a}_p^2)(\sum_{j=0}^{p}b^0_j+\sum_{j=1}^{p-1}{b}_j)-
(a^0_p)^2v_p\right]+{\cO(1)\/z^{2}}\right),\nonumber\\
&\qq \mbox{as}\,\,z\rightarrow\infty.\nonumber
\end{align}

\end{document}